\documentclass[11pt]{amsart}
\usepackage{amssymb,cite,mathrsfs}
\usepackage{esint}
\usepackage[colorlinks=true,urlcolor=blue,
citecolor=red,linkcolor=blue,linktocpage,pdfpagelabels,
bookmarksnumbered,bookmarksopen]{hyperref}
\usepackage[hyperpageref]{backref}
\usepackage[left=2.8cm,right=2.8cm,top=3.2cm,bottom=3.2cm]{geometry}

\numberwithin{equation}{section}

\newtheorem{theorem}{Theorem}[section]
\newtheorem{corollary}[theorem]{Corollary}

\newtheorem{lemma}[theorem]{Lemma}
\newtheorem{proposition}[theorem]{Proposition}

\theoremstyle{definition}

\newtheorem{remark}[theorem]{Remark}

\title[Stability of eigenvalues for the 
fractional $p-$Laplacian]{Stability of variational eigenvalues \\ 
for the fractional $p-$Laplacian}

\author[Brasco]{Lorenzo Brasco}
\author[Parini]{Enea Parini}
\author[Squassina]{Marco Squassina}

\address[L.\ Brasco]{Aix-Marseille Universit\'e, CNRS
\newline\indent
Centrale Marseille, I2M, UMR 7373, 39 Rue Fr\'ed\'eric Joliot Curie
\newline\indent
13453 Marseille, France}
\email{lorenzo.brasco@univ-amu.fr}

\address[E.\ Parini]{Aix-Marseille Universit\'e, CNRS
\newline\indent
Centrale Marseille, I2M, UMR 7373, 39 Rue Fr\'ed\'eric Joliot Curie
\newline\indent
13453 Marseille, France}
\email{enea.parini@univ-amu.fr}

\address[M.\ Squassina]{Dipartimento di Informatica
\newline\indent
Universit\`a degli Studi di Verona
\newline\indent
Verona, Italy}
\email{marco.squassina@univr.it}

\thanks{The research was partially supported by {\it Gruppo Nazionale per l'Analisi Matematica, la Probabilit\`a e le loro Applicazioni} (INdAM) and by the {\it Agence Nationale de la Recherche}, through the project ANR-12-BS01-0014-01 Geometrya.}


\begin{document}
\begin{abstract}
By virtue of $\Gamma-$convergence arguments, we investigate the stability 
of variational eigenvalues associated with a given topological index for the fractional $p-$Laplacian operator, in the
singular limit as the nonlocal operator converges to the $p-$Laplacian. We also obtain the convergence of the corresponding
normalized eigenfunctions in a suitable fractional norm.
\end{abstract}

\maketitle

\begin{center}
\begin{minipage}{11cm}
\small
\tableofcontents
\end{minipage}
\end{center}

\section{Introduction}

\subsection{Overview}
Let $1<p<\infty$, $s \in (0,1)$ and let
$\Omega\subset\mathbb{R}^N$ be a bounded domain with Lipschitz boundary $\partial\Omega$.
Recently, the following nonlocal nonlinear operator was considered in \cite{BraPar,C,FP,IanLiuPerSqu,IaMosSq,meyanez,LL}
\begin{equation}
\label{plap}
(- \Delta_p)^s\, u(x):= 2\, \lim_{\varepsilon \searrow 0} \int_{\mathbb{R}^N \setminus B_\varepsilon(x)}
 \frac{|u(x) - u(y)|^{p-2}\, (u(x) - u(y))}{|x - y|^{N+s\,p}}\, dy, \qquad x \in \mathbb{R}^N.
\end{equation}
For $p=2$, this definition coincides (up to a normalization constant depending on $N$ and $s$, see \cite{CS}) with the linear
fractional Laplacian $(-\Delta)^s$, defined by
\[
(-\Delta)^s=\mathcal{F}^{-1}\,\circ \mathcal{M}_{s}\circ \mathcal{F},
\]
where $\mathcal{F}$ is the Fourier transform operator and $\mathcal{M}_{s}$ is the multiplication by $|\xi|^{2\,s}$.
\par
Many efforts have been devoted to the study of problems involving the fractional $p-$Laplacian operator,
among which we mention eigenvalue problems \cite{BraPar,FP,meyanez,LL}, regularity theory \cite{DKP1,IaMosSq,Kuusi,Kuusi2} and existence of solutions within the framework of Morse theory \cite{IanLiuPerSqu}. For the motivations that lead to
the study of such operators, we refer the reader to the contribution \cite{C} of Caffarelli. For completeness, we also mention that other types of nonlocal quasilinear operators, defined by means of extension properties, can be found in the literature (see \cite{SV}). 
\vskip.2cm\noindent
In this paper, we are concerned with Dirichlet eigenvalues
of $(-\Delta_p)^s$ on the set $\Omega$. These are 
the real (positive) numbers $\lambda$ admitting nontrivial solutions to the following problem
\begin{equation} \label{prob}
\begin{cases}
(- \Delta_p)^s u  =\lambda\, |u|^{p-2}\,u,   & \text{in }\Omega, \\
u  = 0, & \text{in } \mathbb{R}^N \setminus \Omega.
\end{cases}
\end{equation}
It is known that it is possible to construct an infinite sequence of such eigenvalues diverging to $+\infty$. This is done by means of variational methods similar to the so-called {\it Courant minimax principle}, that we briefly recall below. Then our main concern is the study of the singular limit of these {\it variational eigenvalues} as $s\nearrow 1$, in which case the limiting problem of \eqref{prob} is formally given by
\begin{equation} \label{prob-p1}
\begin{cases}
-\Delta_p u=\lambda\, |u|^{p-2}\,u,   & \text{in }\Omega, \\
u  = 0, & \text{on $\partial \Omega$},
\end{cases}
\end{equation}
where $\Delta_p u={\rm div}(|\nabla u|^{p-2}\,\nabla u)$ is the familiar $p-$Laplace operator.
\par
In order to neatly present the subject, we first need some definitions. The natural setting for equations involving the operator $(- \Delta_p)^s$ is the space $W^{s,p}_0(\mathbb{R}^N)$, defined as the completion of $C^\infty_0(\mathbb{R}^N)$ with respect to the standard Gagliardo semi-norm
\begin{equation}
\label{psnorm}
[u]_{W^{s,p}(\mathbb{R}^N)}:=\left(\int_{\mathbb{R}^N}\int_{\mathbb{R}^{N}}\frac{|u(x)-u(y)|^p}{|x-y|^{N+s\,p}}\,dx\,dy\right)^\frac{1}{p}.
\end{equation}
Furthermore, in order to take the Dirichlet condition $u = 0$ 
in $\mathbb{R}^N \setminus \Omega$ into account, we consider the space
$$
\widetilde W_0^{s,p}(\Omega)=\Big\{u:\mathbb{R}^N\to\mathbb{R}\, :\, [u]_{W^{s,p}(\mathbb{R}^N)}<+\infty\,\mbox{ and }\, u=0 \mbox{ in } \mathbb{R}^N \setminus \Omega\Big\},
$$
endowed with \eqref{psnorm}. Since $\Omega$ is Lipschitz, the latter coincides with the space used in \cite{BraLinPar,BraPar} and defined as the completion of $C^\infty_0(\Omega)$ with respect to $[\,\cdot\,]_{W^{s,p}(\mathbb{R}^N)}$ (see Proposition \ref{prop:density} below). Then equation \eqref{prob} has to be intended in the following weak sense:
\[
\int_{\mathbb{R}^N} \int_{\mathbb{R}^N} \frac{|u(x)-u(y)|^{p-2}\, (u(x)-u(y))\,(\varphi(x)-\varphi(y))}{|x-y|^{N+s\,p}}\,dx\,dy=\lambda\, \int_\Omega |u|^{p-2}\,u\,\varphi\,dx,
\]
for every $\varphi\in \widetilde W^{s,p}_0(\Omega)$. Let us introduce
\[
\mathcal{S}_{s,p}(\Omega)=\{u\in \widetilde W^{s,p}_0(\Omega):\|u\|_{L^p(\Omega)}=1\},
\]
and
\[
\mathcal{S}_{1,p}(\Omega)=\{u\in W^{1,p}_0(\Omega):\|u\|_{L^p(\Omega)}=1\},
\]
where $W^{1,p}_0(\Omega)$ is the completion of $C^\infty_0(\Omega)$ with respect to the $L^p$ norm of the gradient.
The {\it $m-$th (variational) eigenvalues} of \eqref{prob} and \eqref{prob-p1} can be obtained as 
\begin{equation}
\label{lambdas}
\lambda^{s}_{m,p}(\Omega):=\inf_{K\in\mathcal W^{s}_{m,p}(\Omega)}\max_{u\in K}\, [u]^p_{W^{s,p}(\mathbb{R}^N)},
\end{equation}
and
\[
\lambda^{1}_{m,p}(\Omega):=\inf_{K\in\mathcal W^{1}_{m,p}(\Omega)}\max_{u\in K}\, \|\nabla u\|^p_{L^p(\Omega)}.
\]
In the previous formulas, we noted for $0<s\le 1$
\begin{equation}
\label{doppiwu}
\mathcal W^{s}_{m,p}(\Omega)=\left\{K\subset\mathcal{S}_{s,p}(\Omega)\,:\,  K \mbox{ symmetric and compact},\, i(K)\ge m\right\},
\end{equation}
and $i(K)$ denotes the {\it Krasnosel'ski\u{\i} genus} of $K$. We recall that for every nonempty and symmetric subset $A\subset X$ of a Banach space, its Krasnosel'ski\u{\i} genus is defined by
\begin{equation}
\label{krasno}
i(A)=\inf\left\{k\in\mathbb{N}\, :\, \exists \mbox{ a continuous odd map } f:A\to\mathbb{S}^{k-1}\right\},
\end{equation}
with the convention that $i(A)=+\infty$, if no such an integer $k$ exists. For completeness, we also mention that for $m=1$ and $m=2$ the previous definitions coincide with
\[
\lambda^s_{1,p}(\Omega)=\min_{u\in \mathcal{S}_{s,p}(\Omega)}[u]^p_{W^{s,p}(\mathbb{R}^N)},\qquad \mbox{\it global minimum},
\]
and
\[
\lambda^s_{2,p}(\Omega)=\inf_{\gamma\in \Sigma(u_1,-u_1)} \max_{u\in \gamma([0,1])}[u]^p_{W^{s,p}(\mathbb{R}^N)},\qquad \mbox{\it mountain pass level},
\]
where $u_1$ is a minimizer associated with $\lambda^s_{1,p}(\Omega)$ and $\Sigma(u_1,-u_1)$ is the set of continuous paths on $\mathcal{S}_{s,p}(\Omega)$ connecting $u_1$ and $-u_1$ (see \cite[Corollary 3.2]{CDG2} for the local case, \cite[Theorem 5.3]{BraPar} for the nonlocal one).
\begin{remark}
For the limit problem \eqref{prob-p1}, the continuity with respect to $p$ of the (variational) eigenvalues $\lambda^{1}_{m,p}$ has been first studied by Lindqvist \cite{stab} and Huang \cite{huang} in the case of the first and second eigenvalue, respectively. Then the problem has been tackled in more generality in \cite{champion_depascale2007, parini2011, littig_schuricht2014}. We also cite the recent paper \cite{DM14} where some generalizations (presence of weights, unbounded sets) have been considered.
\end{remark}

\subsection{Main result}

In order to motivate the investigation pursued in the present paper, it is useful to observe that based upon the results by Bourgain, Brezis and Mironescu \cite{bourgain,BBM2}, we have that if $u\in W_0^{1,p}(\Omega)$ 
\begin{equation}
\label{key}
\lim_{s\nearrow 1} (1-s)\,[u]_{W^{s,p}(\mathbb{R}^N)}^p=K(p,N)\,\|\nabla u\|^p_{L^p(\Omega)},
\end{equation}
(see Proposition~\ref{norm-cont} below).
The constant $K(p,N)$ is given by
\begin{equation}
\label{defK}
K(p,N):=\frac{1}{p}\int_{\mathbb{S}^{N-1}}|\langle \sigma,\mathbf{e}\rangle|^p\, d\mathcal{H}^{N-1}(\sigma),\qquad \mathbf{e}\in \mathbb{S}^{N-1}.
\end{equation}
It is not difficult to see that, due to symmetry reasons, the definition of $K(p,N)$ is indeed independent of the direction $\mathbf{e}\in \mathbb{S}^{N-1}$.
\par
Formula \eqref{key} naturally leads to argue that the nonlocal variational eigenvalues 
$\lambda^{s}_{m,p}$ could converge (once properly renormalized) to the local ones $\lambda^1_{m,p}$. This is the content of the main result of the paper. Observe that we can also assure convergence of the eigenfunctions in suitable (fractional) Sobolev norms.
\begin{theorem}
\label{main}
Let $\Omega\subset\mathbb{R}^N$ be an open and bounded Lipschitz set.
For any $1<p<\infty$ and $m\in\mathbb{N}\setminus\{0\}$ 
\[
\lim_{s\nearrow 1} (1-s)\,  \lambda^s_{m,p}(\Omega)=K(p,N)\,\lambda^1_{m,p}(\Omega).
\]
Moreover, if $u_{s}$ is an eigenfunction of \eqref{prob} corresponding to the variational eigenvalue $\lambda^s_{m,p}(\Omega)$ and such that $\|u_s\|_{L^p(\Omega)}=1$, then there exists a sequence $\{u_{s_k}\}_{k\in\mathbb{N}}\subset\{u_s\}_{s\in(0,1)}$ such that
\[
\lim_{k\to \infty} [u_{s_k}-u]_{W^{t,q}(\mathbb{R}^N)}=0,\qquad \mbox{ for every } p\le q<\infty \mbox{ and every } 0<t<\frac{p}{q},
\]
where $u$ is an eigenfunction of \eqref{prob-p1} corresponding to the variational eigenvalue $\lambda^1_{m,p}(\Omega)$ and such that $\|u\|_{L^p(\Omega)}=1$.
\end{theorem}
\begin{remark}[The case $p=2$]
To the best of our knowledge this result is new already in the linear case $p=2$, namely for the fractional Laplacian operator $(-\Delta)^s$. 
In the theory of stochastic partial differential equations this corresponds
to the case of a stable {\it L\'evy process}. The kernel corresponding to $(-\Delta)^s$ determines the probability distribution of jumps in the value of the stock price, assigning less probability to big jumps as $s$ increases to $1$. Therefore, since the parameter $s$ has to be determined through empirical data, the stability of the spectrum with respect to $s$ allows for more reliable models of random jump-diffusions, see \cite{App} for more details.
\par
It is also useful to recall that for $p=2$, problems \eqref{prob} and \eqref{prob-p1} admit only a discrete set of eigenvalues, whose associated eigenfunctions give an Hilbertian basis of $L^2(\Omega)$ (once properly renormalized). Then we have that these eigenvalues coincide with those defined by \eqref{lambdas}, see Theorem \ref{thm:uguali!} below.
\end{remark}
One of the main ingredients of the proof of Theorem \ref{main} is a $\Gamma-$convergence result for Gagliardo semi-norms, proven in Theorem \ref{thm:gamma} below.  
Namely, by defining the family of functionals $\mathcal E_{s,p}:L^p(\Omega)\to[0,\infty]$ as
\begin{equation}
\label{funzionales}
\mathcal E_{s,p}(u):=\begin{cases} (1-s)^\frac{1}{p}\,[u]_{W^{s,p}(\mathbb{R}^N)},\quad&\mbox{if }u\in \widetilde W^{s,p}_0(\Omega),\\
+\infty\quad&\mbox{otherwise,}
\end{cases}
\end{equation}
and $\mathcal E_{1,p}:L^p(\Omega)\to[0,\infty]$ by 
\begin{equation}
\label{funzionale1}
\mathcal E_{1,p}(u):=\begin{cases}K(p,N)^\frac{1}{p}\,\|\nabla u\|_{L^p(\Omega)},\quad&\mbox{if }u\in W^{1,p}_0(\Omega),\\
+\infty\quad&\mbox{otherwise.}\end{cases}
\end{equation}
we prove that for $s_k\nearrow 1$ we have
\begin{equation}
\label{gliminf}
\mathcal E_{1,p}(u)=\Big(\Gamma-\lim_{k\to\infty}\mathcal E_{s_k,p}\Big)(u),
\qquad\mbox{ for all } u\in L^p(\Omega),
\end{equation}
where $\Gamma-\lim$ denotes the $\Gamma-$limit of functionals, with respect to the norm topology of $L^p(\Omega)$. We refer to \cite{dalmaso1993} for the relevant definitions and facts needed about $\Gamma-$convergence.

\begin{remark}
We point out that a related $\Gamma-$convergence result can be found in the literature, see \cite[Theorem 8]{ponce} by A. Ponce. While his result is for the semi-norms $(1-s)\,[\,\cdot\,]_{W^{s,p}(\Omega)}$ on a bounded set $\Omega$, ours is for the semi-norms $(1-s)\,[\,\cdot\,]_{W^{s,p}(\mathbb{R}^N)}$ on the whole $\mathbb{R}^N$. Moreover, the techniques used in the proofs are slightly different, indeed for the $\Gamma-\liminf$ inequality we follow the one used in \cite{ADePM} for the {\it $s-$perimeter functional}. Such a proof exploits a blow-up technique, introduced by Fonseca and M\"uller in the context of lower-semicontinuity for quasi-convex functionals, see \cite{FM}. As a byproduct of the method, we obtain a variational characterization of the constant $K(p,N)$ appearing in the limit (see Lemma \ref{lm:costanti} below), which is quite typical of the blow-up procedure.
\end{remark}

\begin{remark}
\label{index-r}
In Theorem \ref{main} the variational eigenvalues are defined by means of the Krasnosel'ski\u{\i} genus, but the same result still holds by replacing it
with a general index $i$ having the following properties:
\begin{itemize}
\item[(i)]
$i(K)\in\mathbb{N}\setminus\{0\}$ is defined
whenever $K\not=\emptyset$ is a compact and symmetric subset of a
topological vector space, such that $0\not\in K$;
\item[(ii)]
if $X$ is a topological vector space and
$\emptyset \not =K\subseteq X\setminus\{0\}$ is compact and symmetric,
then there exists $U\subset X\setminus\{0\}$ open set
such that $K\subseteq U$ and
$i(\widehat{K}) \leq i(K)$
for any compact, symmetric and nonempty $\widehat{K}\subseteq U$\,;
\item[(iii)]
if $X, Y$ are two topological vector spaces,
$\emptyset\not=K\subseteq X\setminus\{0\}$ is compact and symmetric
and $\pi:K\to Y\setminus\{0\}$ is continuous and
odd, then $i(\pi(K)) \geq i(K)$\,.
\end{itemize}
Apart from the Krasnosel'ski\u{\i} genus, other examples are the {\it $\mathbb{Z}_2-$cohomological index}
\cite{fadell_rabinowitz1977} and the {\it Ljusternik-Schnirelman Category} \cite[Chapter 2]{St}.
\end{remark}

\subsection{Plan of the paper}
In Section \ref{Prelim}, we collect various preliminary results, such as sharp functional
inequalities and convergence properties in the singular limit $s\nearrow 1$. We point out that even if most of the results of this section are well-known, we need to prove them in order to carefully trace the sharp dependence on the parameter $s$ in all the estimates. 
In Section \ref{sectGamma}, we prove the $\Gamma-$convergence \eqref{gliminf}. For completeness, we also include a convergence result for dual norms, in the spirit of Bourgain-Brezis-Mironescu's result. Then the main result Theorem~\ref{main} is proven in Section \ref{Proof}. Two appendices close the paper and contribute to make it self-contained.

\section{Preliminaries}
\label{Prelim}

\subsection{Some functional inequalities}
We start with an interpolation inequality.
\begin{proposition}[Interpolation inequality]
\label{prop:interpolation}
For every $t\in(0,1)$ and $1<p\le q<r\le+\infty$, we set 
\[
\alpha:=t\,\frac{p}{q}\,\frac{r-q}{r-p}.
\]
Then, for every $u\in C^\infty_0(\mathbb{R}^N)$ and every $0<\beta<\alpha$, we have
\begin{equation}
\label{interpolation}
\begin{split}
\beta^\frac{1}{q}\,[u]_{W^{\beta,q}(\mathbb{R}^N)}&\le C\,\left(\frac{\alpha}{\alpha-\beta}\right)^\frac{1}{q}\,\|u\|_{L^p(\mathbb{R}^N)}^{(1-\theta)\,\left(1-\frac{\beta}{\alpha}\right)}\,\|u\|^\theta_{L^r(\mathbb{R}^N)}\\
&\times \left((1-t)^\frac{1}{p}\,[u]_{W^{t,p}(\mathbb{R}^N)}\right)^{\frac{\beta}{\alpha}\,(1-\theta)}.
\end{split}
\end{equation}
where $C=C(N,p,q)>0$ and $\theta=\theta(p,q,r)\in[0,1)$ is defined by\footnote{For $r=+\infty$, $\alpha$ and $\theta$ are defined accordingly by $\alpha=t\,\frac{p}{q}$ and $\theta=\frac{q-p}{q}$.}
\begin{equation}
\label{theta}
\theta=\frac{r}{r-p}\,\frac{q-p}{q}.
\end{equation}
In the limit case $t=1$, the previous holds in the form 
\begin{equation}
\label{interpolation2}
\begin{split}
\beta^\frac{1}{q}\,[u]_{W^{\beta,q}(\mathbb{R}^N)}&\le C\,\left(\frac{\alpha}{\alpha-\beta}\right)^\frac{1}{q}\,\|u\|_{L^p(\mathbb{R}^N)}^{(1-\theta)\,\left(1-\frac{\beta}{\alpha}\right)}\,\|u\|^\theta_{L^r(\mathbb{R}^N)}\,\|\nabla u\|_{L^p(\mathbb{R}^N)}^{\frac{\beta}{\alpha}\,(1-\theta)},
\end{split}
\end{equation}
\end{proposition}
\begin{proof}
We first consider the case $t\in(0,1)$.
Let $u\in C^\infty_0(\mathbb{R}^N)$, then we have
\[
\begin{split}
[u]^q_{W^{\beta,q}(\mathbb{R}^N)}&=\int_{\{|h|>1\}} \int_{\mathbb{R}^N} \frac{|u(x+h)-u(x)|^q}{|h|^{N+\beta\,q}}\,dx\,dh\\
&+\int_{\{|h|\le 1\}} \int_{\mathbb{R}^N} \frac{|u(x+h)-u(x)|^q}{|h|^{N+\beta\,q}}\,dx\,dh.
\end{split}
\]
The first integral is estimated by
\[
\begin{split}
\int_{\{|h|>1\}} \int_{\mathbb{R}^N} &\frac{|u(x+h)-u(x)|^q}{|h|^{N+\beta\,q}}\,dx\,dh\\
&\le 2^{q-1}\,\int_{\{|h|>1\}}\,\left(\int_{\mathbb{R}^N} \Big(|u(x+h)|^q+|u(x)|^q\Big)\,dx\right)\,\frac{dh}{|h|^{N+\beta\,q}}\\
&=2^{q}\,\int_{\{|h|>1\}}\,\frac{dh}{|h|^{N+\beta\,q}}\, \left(\int_{\mathbb{R}^N} |u|^q\,dx\right)\\
&\le \frac{N\,\omega_N\, 2^q}{\beta\,q}\, \|u\|_{L^p(\mathbb{R}^N)}^{q\,(1-\theta)}\,\|u\|_{L^r(\mathbb{R}^N)}^{q\,\theta},
\end{split}
\]
where $\theta\in[0,1)$ is determined by scale invariance and is given precisely by \eqref{theta}.
In conclusion, 
\begin{equation}
\label{far}
\int_{\{|h|>1\}} \int_{\mathbb{R}^N} \frac{|u(x+h)-u(x)|^q}{|h|^{N+\beta\,q}}\,dx\,dh\le \frac{N\,\omega_N\, 2^q}{\beta\,q}\, \|u\|_{L^p(\mathbb{R}^N)}^{q\,(1-\theta)}\,\|u\|_{L^r(\mathbb{R}^N)}^{q\,\theta}.
\end{equation}
For the other term, for every $\ell>\beta$ we have
\[
\begin{split}
\int_{\{|h|\le 1\}} \int_{\mathbb{R}^N} &\frac{|u(x+h)-u(x)|^q}{|h|^{N+\beta\,q}}\,dx\,dh\\
&=\int_{\{|h|\le 1\}} \left(\int_{\mathbb{R}^N} \frac{|u(x+h)-u(x)|^q}{|h|^{\ell\,q}}\,dx\right)\,\frac{dh}{|h|^{N+(\beta-\ell)\,q}}\\
&\le \int_{\{|h|\le 1\}} \left(\int_{\mathbb{R}^N} \frac{|u(x+h)-u(x)|^p}{|h|^{\frac{\ell}{1-\theta}\,p}}\,dx\right)^{q\,\frac{1-\theta}{p}}\\
&\times\left(\int_{\mathbb{R}^N} |u(x+h)-u(x)|^{r}\,dx\right)^{q\,\frac{\theta}{r}}\,\frac{dh}{|h|^{N+(\beta-\ell)\,q}}.
\end{split}
\]
We choose 
\begin{equation}
\label{elle}
\frac{\ell}{1-\theta}=t,\qquad \mbox{ i.e. }\quad \ell=t\,\frac{p}{q}\,\frac{r-q}{r-p}=\alpha,
\end{equation}
and use \cite[Lemma A.1]{BraLinPar}, i.e. 
\begin{equation}
\label{quozienti}
\int_{\mathbb{R}^N} \frac{|u(x+h)-u(x)|^p}{|h|^{t\,p}}\,dx\le C\,(1-t)\,[u]^p_{W^{t,p}(\mathbb{R}^N)},
\end{equation}
for some $C=C(N,p)>0$.
On the other hand, we have
\[
\left(\int_{\mathbb{R}^N} |u(x+h)-u(x)|^{r}\,dx\right)^{q\,\frac{\theta}{r}}\le 2^{q\,\theta}\,\left(\int_{\mathbb{R}^N} |u|^{r}\,dx\right)^{q\,\frac{\theta}{r}}.
\]
Thus we get
\begin{equation}
\label{near}
\begin{split}
\int_{\{|h|\le 1\}} \int_{\mathbb{R}^N} &\frac{|u(x+h)-u(x)|^q}{|h|^{N+\beta\,q}}\,dx\,dh\\
&\le \frac{C}{\alpha-\beta}\,2^{q\,\theta}\,\|u\|_{L^{r}(\mathbb{R}^N)}^{q\,\theta}\,\left((1-t)^\frac{1}{p}\, [u]_{W^{t,p}(\mathbb{R}^N)}\right)^{q\,(1-\theta)},
\end{split}
\end{equation}
where $C=C(N,p,q)>0$. 
By combining \eqref{far} and \eqref{near}, we get
\[
\begin{split}
[u]^q_{W^{\beta,q}(\mathbb{R}^N)}&\le \frac{C}{\beta}\,\|u\|_{L^p(\mathbb{R}^N)}^{q\,(1-\theta)}\,\|u\|_{L^r(\mathbb{R}^N)}^{q\,\theta}\\
&+\frac{C}{\alpha-\beta}\,\|u\|_{L^{r}(\mathbb{R}^N)}^{q\,\theta}\,\left((1-t)^\frac{1}{p}\, [u]_{W^{t,p}(\mathbb{R}^N)}\right)^{q\,(1-\theta)},
\end{split}
\]
possibly with a different $C=C(N,p,q)>0$.
We now use the previous inequality with $u_\chi(x)=u(x\,\chi^{-1/q})$ and optimize in $\chi>0$. We get
\begin{equation}
\label{bordello}
\begin{split}
\chi^{\alpha-\beta}\,[u]^q_{W^{\beta,q}(\mathbb{R}^N)}&- \frac{C\,\chi^\alpha}{s}\,\|u\|_{L^p(\mathbb{R}^N)}^{q\,(1-\theta)}\,\|u\|_{L^r(\mathbb{R}^N)}^{q\,\theta}\\
&\le \frac{C}{\alpha-\beta}\,\|u\|_{L^{r}(\mathbb{R}^N)}^{q\,\theta}\,\left((1-t)^\frac{1}{p}\, [u]_{W^{t,p}(\mathbb{R}^N)}\right)^{q\,(1-\theta)},
\end{split}
\end{equation}
still for some constant $C=C(N,p,q)>0$. Observe that by hypothesis on $s$ we have $\alpha-\beta>0$. The left-hand side of \eqref{bordello} is maximal for
\[
\chi_0=\left(\frac{\alpha-\beta}{\alpha}\,\frac{\beta}{C}\,\frac{[u]^q_{W^{\beta,q}(\mathbb{R}^N)}}{\|u\|_{L^p(\mathbb{R}^N)}^{q\,(1-\theta)}\,\|u\|_{L^r(\mathbb{R}^N)}^{q\,\theta}} \right)^\frac{1}{\beta}.
\]
Thus we get
\[
\begin{split}
\left(\frac{\alpha-\beta}{\alpha}\,\frac{\beta}{C}\right)^\frac{\alpha-\beta}{\beta}\, \frac{\beta}{\alpha}\, \frac{[u]^{q\,\frac{\alpha}{\beta}}_{W^{\beta,q}(\mathbb{R}^N)}}{\|u\|_{L^p(\mathbb{R}^N)}^{q\,(1-\theta)\,\frac{\alpha-\beta}{\beta}}\,\|u\|_{L^{r}(\mathbb{R}^N)}^{q\,\theta\,\frac{\alpha-\beta}{\beta}}}&\le \frac{C}{\alpha-\beta}\,\|u\|_{L^{r}(\mathbb{R}^N)}^{q\,\theta}\\
&\times\left((1-t)^\frac{1}{p}\, [u]_{W^{t,p}(\mathbb{R}^N)}\right)^{q\,(1-\theta)},
\end{split}
\]
that is
\[
\begin{split}
[u]_{W^{s,q}(\mathbb{R}^N)}\le \left(\frac{C}{\beta}\, \frac{\alpha}{\alpha-\beta}\right)^\frac{1}{q}\|u\|_{L^p(\mathbb{R}^N)}^{(1-\theta)\left(1-\frac{\beta}{\alpha}\right)}\|u\|^\theta_{L^r(\mathbb{R}^N)}\,\left((1-t)^\frac{1}{p}[u]_{W^{t,p}(\mathbb{R}^N)}\right)^{\frac{\beta}{\alpha}\,(1-\theta)},
\end{split}
\]
with $C=C(N,p,q)>0$.
\vskip.2cm\noindent
In order to prove \eqref{interpolation2}, it is sufficient to repeat the previous proof, this time replacing the choice \eqref{elle} by
\[
\ell=\frac{p}{q}\, \frac{r-q}{r-p},\qquad \mbox{ so that }\quad \frac{\ell}{1-\theta}=1,
\]
and then using that
\[
\int_{\mathbb{R}^N} \frac{|u(x+h)-u(x)|^p}{|h|^{p}}\,dx\le \int_{\mathbb{R}^N} |\nabla u|^p\,dx,
\]
in place of \eqref{quozienti}, which follows from basic calculus and invariance by translations of the $L^p$ norm.
\end{proof}
\begin{corollary}
\label{cor:inclusione}
Let $1<p<\infty$ and $s\in(0,1)$. Then, for every $u\in C^\infty_0(\mathbb{R}^N)$, we have
\begin{equation}
\label{interpolazione_weak}
s\,(1-s)\,[u]^p_{W^{s,p}(\mathbb{R}^N)}\le C\, \|u\|_{L^p(\mathbb{R}^N)}^{(1-s)\,p}\,\|\nabla u\|_{L^p(\mathbb{R}^N)}^{s\,p},
\end{equation}
for some constant $C=C(N,p)>0$.
In particular, if  $\Omega\subset\mathbb{R}^N$ is an open bounded set with Lipschitz boundary, then we have $W^{1,p}_0(\Omega)\subset W^{s,p}_0(\mathbb{R}^N)$ and
\begin{equation}
\label{immersione}
s\,(1-s)\,[u]^p_{W^{s,p}(\mathbb{R}^N)}\le C\,\left(\lambda^{1}_{1,p}(\Omega)\right)^{s-1}\,\|\nabla u\|^p_{L^p(\Omega)},\qquad u\in W^{1,p}_0(\Omega),
\end{equation}
where $C=C(N,p)>0$. 
\end{corollary}
\begin{proof}
In order to prove \eqref{interpolazione_weak}, it is sufficient to use \eqref{interpolation2} with $t=1$, $\beta=s$ and $q=p$: observe that in this case $\alpha=1$ and $\theta=0$.
\vskip.2cm\noindent
For $u\in C^\infty_0(\Omega)$, by inequality \eqref{interpolazione_weak} we get 
\[
[u]^p_{W^{s,p}(\mathbb{R}^N)}\le \frac{C}{s\,(1-s)}\,\|u\|^{(1-s)\,p}_{L^p(\Omega)}\,\|\nabla u\|^{s\,p}_{L^p(\Omega)}.
\]
If we now apply the Poincar\'e inequality for $W^{1,p}_0(\Omega)$ on the right-hand side,
we obtain inequality \eqref{immersione} for functions in $C^\infty_0(\Omega)$. By using density of $C^\infty_0(\Omega)$ in the space $W^{1,p}_0(\Omega)$, we get the desired conclusion. \end{proof}
In what follows, we define the sharp Sobolev constant
\begin{equation}
\label{talentiana}
T_{p,s}:=\sup_{u\in W^{s,p}_0(\mathbb{R}^N)} \left\{\left(\int_{\mathbb{R}^N} |u|^\frac{N\,p}{N-s\,p}\,dx\right)^\frac{N-s\,p}{N}\, :\, [u]_{W^{s,p}(\mathbb{R}^N)}=1\right\}.
\end{equation}
We need the following result by Maz'ya and Shaponishkova, see \cite[Theorem 1]{MazSha}. \begin{theorem}[Estimate of the sharp constant]
\label{thm:mazza}
Let $1<p<\infty$ and $s\in(0,1)$ be such that $s\,p<N$. 
Then for the constant \eqref{talentiana} we have the estimate
\[
T_{p,s}\le \mathcal{T}\,\frac{s\,(1-s)}{(N-s\,p)^{p-1}},
\]
for some $\mathcal{T}=\mathcal{T}(N,p)>0$. 
\end{theorem}

\begin{theorem}[Hardy inequality for convex sets]
\label{Hardy-Opt}
Let $1<p<\infty$ and $s\in (0,1)$ with $s\,p>1$. Then for any convex domain $\Omega\subset\mathbb{R}^N$ and every  $u\in C^\infty_0(\Omega)$ we have
\begin{equation}
\label{hardysubopt}
\left(\frac{s\,p-1}{p}\right)^p\,\mathcal{C}_{N,p}\,\left\|\frac{u}{\delta_\Omega^{s}}\right\|^p_{L^p(\Omega)}\le (1-s)\,[u]^p_{W^{s,p}(\Omega)}.
\end{equation}
with $\delta_\Omega(x)={\rm dist}(x, \partial \Omega)$ and ${\mathcal C}_{N,p}>0$ a costant depending on $N$ and $p$ only. 
\end{theorem}
\begin{proof}
The Hardy inequality without sharp constant can be found in \cite[Theorem 1.1]{Dy}. 
The sharp constant for convex sets was obtained in \cite[Theorem 1.2]{LossSloane}, where it is proven
\[
{\mathcal D}_{N,p,s}\,\left\|\frac{u}{\delta_\Omega^{s}}\right\|^p_{L^p(\Omega)}\le (1-s)\,[u]^p_{W^{s,p}(\Omega)}.
\]
The optimal constant ${\mathcal D}_{N,p,s}$ is given by
$$
{\mathcal D}_{N,p,s}=2\,\pi^{\frac{N-1}{2}}\,\frac{\Gamma(\frac{1+s\,p}{2})}{\Gamma(\frac{N+s\,p}{2})}
\int_0^1\frac{\left(1-r^{\frac{s\,p-1}{p}} \right)^p}{(1-r)^{1+s\,p}}\,dr.
$$
We claim that
\begin{equation}
\label{stima}
{\mathcal D}_{N,p,s}\ge \left(\frac{s\,p-1}{p}\right)^p\,\frac{\mathcal{C}_{N,p}}{1-s}.
\end{equation}
Indeed, by concavity of the map $\tau\mapsto \tau^{(s\,p-1)/p}$ we have
\[
1-r^\frac{s\,p-1}{p}\ge \frac{s\,p-1}{p}\,(1-r),\qquad r\ge 0.
\]
Thus we get
\[
\int_0^1\frac{\left(1-r^{\frac{s\,p-1}{p}} \right)^p}{(1-r)^{1+s\,p}}\,dr\ge \left(\frac{s\,p-1}{p}\right)^p\, \int_0^1 (1-r)^{p-1-s\,p}\,dr=\left(\frac{s\,p-1}{p}\right)^p\,\frac{1}{p\,(1-s)}.
\]
On the other hand, from the definition of $\Gamma$ we also have
\[
\Gamma\left(\frac{1+s\,p}{2}\right)\ge \int_0^1 t^\frac{s\,p-1}{2}\, e^{-t}\,dt\ge \int_0^1 t^\frac{p-1}{2}\, e^{-t}\,dt=:c_1(p)>0,
\]
and also
\[
\Gamma\left(\frac{N+s\,p}{2}\right)\le \int_0^1 t^\frac{N-1}{2}\, e^{-t}\,dt+\int_1^{+\infty} t^\frac{N-2+p}{2}\, e^{-t}\,dt=:c_2(p)>0.
\]
By using these estimates, we get \eqref{stima}.
\end{proof}
The next result is a Poincar\'e inequality for Gagliardo semi-norms. This is classical, but as always we have to carefully trace the sharp dependence on $s$ of the constants concerned.
\begin{proposition}[Poincar\'e inequalities]
\label{lm:spoin}
Let $s\in(0,1)$ and $1\le p<\infty$. Let $\Omega\subset\mathbb{R}^N$ be an open and bounded set. Then
\begin{equation}
\label{spoin}
\|u\|^p_{L^p(\Omega)}\le C\,\mathrm{diam}(\Omega)^{s\,p}\, (1-s)\, [u]_{W^{s,p}(\mathbb{R}^N)}^p,\qquad u\in C^\infty_0(\Omega),
\end{equation}
for a constant $C=C(N,p)>0$.
Moreover, if $\Omega$ is convex and $s\,p>1$, then 
\begin{equation}
\label{spoinbis}
\|u\|^p_{L^p(\Omega)}\le \frac{C}{(s\,p-1)^p}\,\mathrm{diam}(\Omega)^{s\,p}\, (1-s)\, [u]_{W^{s,p}(\Omega)}^p,\qquad u\in C^\infty_0(\Omega),
\end{equation}
possibly with a different constant $C=C(N,p)>0$, still independent of $s$.
\end{proposition}
\begin{proof}
Since $\Omega$ is bounded, we have $\Omega\subset B_R(x_0)$, with $2\,R=\mathrm{diam}(\Omega)$ and $x_0\in\Omega$. Let $h\in\mathbb{R}^N$ be such that $|h|>2\,R$, so that
\[
x+h\in \mathbb{R}^N\setminus \Omega,\qquad \mbox{ for every }x\in \Omega.
\]
Then for every $u\in C^\infty_0(\Omega)$ we have
\[
\begin{split}
\int_\Omega |u(x)|^p\, dx=\int_\Omega |u(x+h)-u(x)|^p\, dx&\le \int_{\mathbb{R}^N} |u(x+h)-u(x)|^p\, dx\\
&=|h|^{s\,p}\,\int_{\mathbb{R}^N} \frac{|u(x+h)-u(x)|^p}{|h|^{s\,p}}\, dx\\
&\le |h|^{s\,p} C_{N,p}\, (1-s)\, [u]^p_{W^{s,p}(\mathbb{R}^N)},
\end{split}
\]
where in the last estimate we used \cite[Lemma A.1]{BraLinPar}. By taking the infimum on the admissible $h$, we get \eqref{spoin}.
\par
Let us now suppose that $\Omega$ is convex and $s\,p>1$. In order to prove \eqref{spoinbis}, we proceed as in the proof of \cite[Proposition B.1]{BraLinPar}. For every $u\in C^\infty_0(\Omega)$ we have
\[
[u]^p_{W^{s,p}(\mathbb{R}^N)}=[u]^p_{W^{s,p}(\Omega)}+2\, \int_\Omega \int_{\mathbb{R}^N\setminus\Omega} \frac{|u(x)|^p}{|x-y|^{N+s\,p}}\,dx\,dy.
\]
In order to estimate the last term, we first observe that, if $\delta_\Omega(x)={\rm dist}(x, \partial \Omega)$, we get 
\[
\begin{split}
\int_\Omega \int_{\mathbb{R}^N\setminus\Omega} \frac{|u(x)|^p}{|x-y|^{N+s\,p}}\,dx\,dy&\le \int_\Omega \left(\int_{\mathbb{R}^N\setminus B_{\delta_\Omega(x)}(x)}\frac{1}{|x-y|^{N+s\,p}}\,dy\right)\,|u(x)|^p\,dx\\
&= N\,\omega_N\,\int_\Omega \left(\int_{\delta_\Omega(x)}^{+\infty} \varrho^{-1-s\,p}\,d\varrho\right)\,|u(x)|^p\,dx\\
&\le N\,\omega_N\, \int_\Omega \frac{|u|^p}{\delta_\Omega^{s\,p}}\,dx,
\end{split}
\]
where we also used that $s\,p>1$.
We can now use Hardy inequality \eqref{hardysubopt}, so to obtain
\[
\int_\Omega \int_{\mathbb{R}^N\setminus\Omega} \frac{|u(x)|^p}{|x-y|^{N+s\,p}}\,dx\,dy\le N\,\omega_N\,\left(\frac{p}{s\,p-1}\right)^p\,\frac{1-s}{\mathcal{C}_{N,p}}\,[u]^p_{W^{s,p}(\Omega)}.
\]
By also using that $1-s<1$, we finally get
\[
[u]^p_{W^{s,p}(\mathbb{R}^N)}\le \left[1+2\,N\,\omega_N\,\left(\frac{p}{s\,p-1}\right)^p\,\frac{1}{\mathcal{C}_{N,p}}\right]\,[u]^p_{W^{s,p}(\Omega)}.
\]
By combining this, \eqref{spoin} and observing that $p/(s\,p-1)>1$, we get \eqref{spoinbis}.
\end{proof}
\begin{remark}
The constant in \eqref{spoinbis} degenerates as $s\,p$ goes to $1$. Indeed, for $s\,p\le 1$ a Poincar\'e inequality like \eqref{spoinbis} is not possible (see \cite[Remark 2.7]{BraPar}).
\end{remark}

\subsection{From nonlocal to local}
We will systematically use the following result.
\begin{theorem}[\cite{bourgain}]
\label{teo:BBM}
Let $1<p<\infty$ and let $\Omega\subset\mathbb{R}^N$ be an open bounded set with Lipschitz boundary. Then, for every $u\in W^{1,p}(\Omega)$, we have 
\begin{equation}
\label{convBBM}
\lim_{s\nearrow 1} (1-s)\,[u]_{W^{s,p}(\Omega)}^p=K(p,N)\,\|\nabla u\|^p_{L^p(\Omega)},
\end{equation}
where $K(p,N)$ is defined in \eqref{defK}. 
\end{theorem}
More precisely, we will need the following extension. The main difference with Theorem \ref{teo:BBM} is that functions are now taken in $W^{1,p}_0(\Omega)$ and the seminorm on $\Omega$ is replaced by the seminorm on the whole $\mathbb{R}^N$.
\begin{proposition}
\label{norm-cont}
Let $1<p<\infty$ and let $\Omega\subset\mathbb{R}^N$ be an open bounded set with Lipschitz boundary. Assume that $u\in W^{1,p}_0(\Omega)$, then we have 
\begin{equation}
\label{conv-s}
\lim_{s\nearrow 1} (1-s)\,[u]^p_{W^{s,p}(\mathbb{R}^N)}=K(p,N)\,\|\nabla u\|^p_{L^p(\Omega)},
\end{equation}
where $K(p,N)>0$ is defined in \eqref{defK}.
\end{proposition}
\begin{proof}
Let $u\in W^{1,p}_0(\Omega)$, we observe that $u\in  W^{s,p}_0(\mathbb{R}^N)$ for all $s\in (0,1)$ thanks to Corollary \ref{cor:inclusione}. Furthermore, since $\Omega$ is a bounded set, by virtue of Theorem \ref{teo:BBM} we have
\[
\lim_{s\nearrow 1} (1-s)\, [u]^p_{W^{s,p}(\Omega)}
=K(p,N)\,\|\nabla u\|_{L^p(\Omega)}.
\]
Let us first prove that \eqref{conv-s} holds for $u\in C^\infty_0(\Omega)$.
Recalling that $u=0$ outside $\Omega$, we have
\[
(1-s)\,[u]^p_{W^{s,p}(\mathbb{R}^N)}=
(1-s)\,[u]^p_{W^{s,p}(\Omega)}+2\,(1-s)\int_{\Omega}\int_{\mathbb{R}^{N}\setminus\Omega}\frac{|u(x)|^p}{|x-y|^{N+s\,p}}\,dx\,dy.
\]
This yields
\[
\begin{split}
\lim_{s\nearrow 1}(1-s)\,[u]^p_{W^{s,p}(\mathbb{R}^N)}&=
K(p,N)\,\int_{\Omega}|\nabla u|^p\,dx\\
&+2\lim_{s\nearrow 1}(1-s)\,\int_{\Omega}
\int_{\mathbb{R}^N\setminus\Omega}\frac{|u(x)|^p}{|x-y|^{N+s\,p}}\,dy\,dx.
\end{split}
\]
Since $u\in C^\infty_0(\Omega)$, we have ${\rm dist}(\partial K,\partial\Omega)>0$
where $K$ is the support of $u$. It follows that
$$
\int_{\Omega}
\left(\int_{\mathbb{R}^{N}\setminus\Omega}\frac{|u(x)|^p}{|x-y|^{N+sp}}dy\right)dx\leq
\|u\|_{L^p(\Omega)}^p\,\int_{\mathbb{R}^{N}\setminus\Omega}\frac{1}{\delta_K(y)^{N+s\,p}}\,dy,
$$
where we set $\delta_K(y)=\mathrm{dist\,}(y,\partial K)$.
Hence, there exists a constant $C=C(N,p)>0$, such that if $R=\mathrm{dist}(\mathbb{R}^N\setminus\Omega,\partial K)>0$, then we have
$$
\int_{\mathbb{R}^{N}\setminus\Omega}\frac{1}{\delta_K(y)^{N+s\,p}}\,dy\leq
\int_{\mathbb{R}^{N}\setminus B(0,R)}\frac{1}{|y|^{N+s\,p}}\,dy=\frac{C}{s}\,R^{-s\,p}.
$$
It follows that
$$
\lim_{s\nearrow 1}(1-s)\int_{\Omega}
\left(\int_{\mathbb{R}^{N}\setminus\Omega}\frac{|u(x)|^p}{|x-y|^{N+sp}}dy\right)dx=0,
$$
and the claim is proved for $u\in C^\infty_0(\Omega)$.
\vskip.2cm\noindent
Assume now that $u\in W^{1,p}_0(\Omega)$. Then there exists a sequence $\{\phi_j\}_{j\in\mathbb{N}}\subset C^\infty_0(\Omega)$ such that $\|\nabla \phi_j-\nabla u\|_{L^p(\Omega)}
\to 0$ as $j$ goes to $\infty$. In turn,
by inequality \eqref{immersione} we have 
$$
(1-s)^\frac{1}{p}\,[\phi_j-u]_{W^{s,p}(\mathbb{R}^N)}\leq C\, \|\nabla\phi_j-\nabla u\|_{L^p(\Omega)},
$$
with $C$ independent of $s$ and $j$.
Thus for every $\varepsilon>0$, there exists $j_0\in\mathbb{N}$ independent of $s$ such that
\begin{equation}
\label{tardi}
\Big|\|\nabla u\|_{L^p(\Omega)}-\|\nabla \phi_j\|_{L^p(\Omega)}\Big|\le \|\nabla \phi_j-\nabla u\|_{L^p(\Omega)}\le \varepsilon,
\end{equation}
and consequently
\[
\left|(1-s)^\frac{1}{p}\,[\phi_j]_{W^{s,p}(\mathbb{R}^N)}-(1-s)^\frac{1}{p}\,[u]_{W^{s,p}(\mathbb{R}^N)}\right|\le (1-s)^\frac{1}{p}\,[\phi_j-u]_{W^{s,p}(\mathbb{R}^N)}\le C\,\varepsilon,
\]
for every $j\geq j_0$. Then for every $j\ge j_0$
\[
\begin{split}
(1-s)^\frac{1}{p}\,[\phi_j]_{W^{s,p}(\mathbb{R}^N)}-C\,\varepsilon&\le (1-s)^\frac{1}{p}\,[u]_{W^{s,p}(\mathbb{R}^N)} \leq (1-s)^\frac{1}{p}\,[\phi_j]_{W^{s,p}(\mathbb{R}^N)}+C\,\varepsilon,
\end{split}
\]
for every $j\geq j_0$. By using the first part of the proof we thus get for every $j\ge j_0$
\[
\begin{split}
K(p,N)^\frac{1}{p}\, \|\nabla \phi_j\|_{L^p(\Omega)}-C\,\varepsilon&\le \lim_{s\nearrow 1} (1-s)^\frac{1}{p}\, [u]_{W^{s,p}(\mathbb{R}^N)}\\
&\le K(p,N)^\frac{1}{p}\, \|\nabla \phi_j\|_{L^p(\Omega)}+C\,\varepsilon.
\end{split}
\]
If we now use \eqref{tardi} and exploit the arbitrariness of $\varepsilon>0$, we get \eqref{conv-s} for a general $u\in W^{1,p}_0(\Omega)$.
\end{proof}

\subsection{Dual spaces}

Let $\Omega\subset\mathbb{R}^N$ be as always an open and bounded set with Lipschitz boundary. Let $1<p<\infty$ and $s\in(0,1)$, we set $p'=p/(p-1)$ and 
\[
W^{-s,p'}(\Omega)=\left\{F:\widetilde W^{s,p}_0(\Omega)\to\mathbb{R}\, :\, F \mbox{ linear and continuous}\right\},
\]
which is equipped with the natural dual norm
\[
\|F\|_{W^{-s,p'}(\Omega)}=\sup_{u\in \widetilde W^{s,p}_0(\Omega)\setminus\{0\}} \frac{|\langle F,u\rangle|}{[u]_{W^{s,p}(\mathbb{R}^N)}}.
\]
The symbol $\langle\cdot,\cdot\rangle$ denotes the relevant duality product. For $s=1$, the space $W^{-1,p'}(\Omega)$ and the corresponding dual norm are defined accordingly.
\par
The following is the dual version of \eqref{immersione}.
\begin{lemma}
Let $F\in W^{-s,p'}(\Omega)$, then we have
\begin{equation}
\label{mandorlini}
\|F\|_{W^{-1,p'}(\Omega)}\le \left(\frac{C}{s\,(1-s)}\right)^\frac{1}{p}\,\left(\lambda^1_{1,p}(\Omega)\right)^\frac{s-1}{p}\, \|F\|_{W^{-s,p'}(\Omega)}.
\end{equation}
\end{lemma}
\begin{proof}
By \eqref{immersione} we have 
\[
[u]_{W^{s,p}(\mathbb{R}^N)}\le \left(\frac{C}{s\,(1-s)}\right)^\frac{1}{p}\,\left(\lambda^{1}_{1,p}(\Omega)\right)^\frac{s-1}{p}\,\|\nabla u\|_{L^p(\Omega)},\qquad u\in C^\infty_0(\Omega).
\]
Thus for every $u\in C^\infty_0(\Omega)\setminus\{0\}$ we have
\[
\frac{|\langle F,u\rangle|}{[u]_{W^{s,p}(\mathbb{R}^N)}}\ge \left(\frac{s\,(1-s)}{C}\right)^\frac{1}{p}\,\left(\lambda^{1}_{1,p}(\Omega)\right)^\frac{1-s}{p}\, \frac{|\langle F,u\rangle|}{\|\nabla u\|_{L^p(\Omega)}}.
\]
By taking the supremum over $u$, the conclusion follows from the definition of dual norm.
\end{proof}
If $F\in W^{-s,p'}(\Omega)$, by a simple homogeneity argument (i.e. replacing $u$ by $t\,u$ and then optimizing in $t$) we have
\[
\begin{split}
\max_{u\in \widetilde W^{s,p}_0(\Omega)} \left\{\langle F,u\rangle-(1-s)\,[u]^p_{ W^{s,p}(\mathbb{R}^N)}\right\}
&=\left(\frac{1}{p}\right)^\frac{1}{p-1}\,\frac{1}{p'}\, \left(\frac{\|F\|_{W^{-s,p'}(\Omega)}}{(1-s)^\frac{1}{p}}\right)^{p'}.
\end{split}
\]
Thus in particular we get
\begin{equation}
\label{duale}
\frac{\|F\|_{W^{-s,p'}(\Omega)}}{(1-s)^\frac{1}{p}}=p^\frac{1}{p}\,\left(-p'\,\min_{u\in \widetilde W^{s,p}_0(\Omega)} \left\{(1-s)\,[u]^p_{ W^{s,p}(\mathbb{R}^N)}-\langle F,u\rangle\right\}\right)^\frac{1}{p'}.
\end{equation}
In the local case $s=1$, with similar computations we get
\begin{equation}
\label{duale1}
\frac{\|F\|_{W^{-1,p'}(\Omega)}}{K(p,N)^\frac{1}{p}}=p^\frac{1}{p}\,\left(-p'\,\min_{u\in \widetilde W^{s,p}_0(\Omega)} \left\{K(p,N)\,\|\nabla u\|^p_{L^p(\Omega)}-\langle F,u\rangle\right\}\right)^\frac{1}{p'}.
\end{equation}

\subsection{A bit of regularity} 
We conclude this section with a regularity result. This is not new, but once again our main concern is the dependence on $s$ of the constants entering in the estimates below. We also need to pay particular attention to the case $p=N$, which is borderline in the limit as $s$ goes $1$.
\begin{theorem}
Let $1<p<\infty$ and $s\in(0,1)$. If $u\ \in \widetilde W^{s,p}_0(\Omega)$ is an eigenfunction of $(-\Delta_p)^s$ with eigenvalue $\lambda$, then we have:
\begin{itemize}
\item if $1<p<N$
\begin{equation}
\label{boundep}
\|u\|_{L^\infty(\Omega)}\le \left[\frac{C}{(N-s\,p)^{p-1}}\,s\,(1-s)\,\lambda\right]^\frac{N}{s\,p^2}\,\|u\|_{L^p(\Omega)},
\end{equation}
for a constant $C=C(N,p)>0$;
\vskip.2cm
\item if $p=N$ and $s\ge 3/4$
\begin{equation}
\label{cresima}
\|u\|_{L^\infty}\le \left[C\,\Big(\mathrm{diam}(\Omega)\Big)^{N\,\frac{2\,s-1}{2}}\,s\,(1-s)\,\lambda\right]^\frac{2}{N}\,\|u\|_{L^N(\Omega)},
\end{equation}
for a constant $C=C(N)>0$;
\vskip.2cm
\item if $p>N$ and $s\in(0,1)$ is such that $s\,p>N$
\begin{equation}
\label{bounder}
\|u\|_{L^\infty(\Omega)}\le \Big[C\,(1-s)\,\lambda\Big]^\frac{1}{p}\,\|u\|_{L^p(\Omega)}\,\mathrm{diam\,}(\Omega)^{s-\frac{N}{p}},
\end{equation}
for a constant $C=C(N,p)>0$.
\end{itemize}
\end{theorem}
\begin{proof}
In the case $1<p<N$, we have of course $s\,p<N$ as well. Then by appealing to \cite[Theorem 3.3 \& Remark 3.4]{BraLinPar} and \cite[Remark 3.2]{BraPar} we know that
\[
\|u\|_{L^\infty(\Omega)}\le \left[\left(\frac{N}{N-s\,p}\right)^{\frac{N-s\,p}{s}\,\frac{p-1}{p}}\,T_{p,s}\,\lambda\right]^\frac{N}{s\,p^2}\,\|u\|_{L^p(\Omega)}.
\]
where $T_{p,s}$ is the sharp Sobolev constant \eqref{talentiana}. By using Theorem \ref{thm:mazza}, we obtain
\[
\|u\|_{L^\infty(\Omega)}\le \left[\left(\frac{1}{1-(s\,p)/N}\right)^{\left(1-\frac{s\,p}{N}\right)\frac{N}{s\,p}\,(p-1)}\,\frac{\mathcal{T}\,s}{(N-s\,p)^{p-1}}\,(1-s)\,\lambda\right]^\frac{N}{s\,p^2}\,\|u\|_{L^p(\Omega)}.
\]
Then from the previous we get \eqref{boundep}, once it is noticed that $(1-(s\,p)/N)^{1-(s\,p)/N}$ is bounded from below by a universal constant $c>0$, for $0<s\,p<N$.
\vskip.2cm\noindent
Let us now consider the case 
\[
p=N\qquad \mbox{ and }\qquad \frac{3}{4}\le s.
\] 
By proceeding as in the second part of the proof of\footnote{This is based on the Moser's iteration technique, this time with the Sobolev embedding
\[
\widetilde W^{s,N}_0(\Omega)\hookrightarrow L^{2\,N}(\Omega)\qquad \mbox{ in place of }\qquad \widetilde W^{s,N}_0(\Omega)\hookrightarrow L^\frac{N}{1-s}(\Omega).
\] 
The proof in \cite[Theorem 3.3]{BraLinPar} is for the first eigenfunction, but it can be easily adapted to the case of any eigenfunction, as observed in \cite[Remark 3.2]{BraPar}.} \cite[Theorem 3.3]{BraLinPar} and using the same notation, after a Moser's iteration we end up with
\begin{equation}
\label{alpelo}
\|u\|_{L^\infty}\le C\, \left(\frac{\lambda}{\alpha^s_N(\Omega)}\right)^\frac{2}{N}\,\|u\|_{L^N(\Omega)},
\end{equation}
where $C=C(N)>0$ and the geometric constant $\alpha^s_N(\Omega)$ is given by
\[
\alpha^s_N(\Omega)=\min_{u\in\widetilde W^{s,N}_0(\Omega)} \Big\{[u]^N_{W^{s,N}(\mathbb{R}^N)}\, :\, \|u\|_{L^{2\,N}(\Omega)}=1\Big\}.
\]
Observe that this is not $0$, since $\widetilde W^{s,N}_0(\Omega)\hookrightarrow L^{2\,N}(\Omega)$ as soon as 
\[
2\,N<\frac{N}{1-s}\qquad \Longleftrightarrow \qquad \frac{1}{2}<s,
\]
which is verified. In order to estimate $\alpha^s_N(\Omega)$ from below, we observe that by choosing $q=p=N$, $t=s$ and $\beta=1/2$ in \eqref{interpolation},
we get ($C$ denotes a constant depending on $N$ only, varying from a line to another)
\[
\begin{split}
\frac{1}{2}\,[u]^N_{W^{\frac{1}{2},N}(\mathbb{R}^N)}&\le C\,\left(\dfrac{2\,s}{2\,s-1}\right)\,\|u\|_{L^N(\mathbb{R}^N)}^{N\,\left(1-\frac{1}{2\,s}\right)}
\,\left((1-s)\,[u]^N_{W^{s,N}(\mathbb{R}^N)}\right)^{\frac{1}{2\,s}}\\
&\le C\,\left(\dfrac{2\,s}{2\,s-1}\right)\,\Big(\mathrm{diam}(\Omega)\Big)^{N\,\frac{2\,s-1}{2}}
\,\left((1-s)\,[u]^N_{W^{s,N}(\mathbb{R}^N)}\right),
\end{split}
\]
where in the second inequality we used Poincar\'e inequality \eqref{spoin}. We can now use Sobolev inequality in the left-hand side, i.e.
\[
[u]^N_{W^{\frac{1}{2},N}(\mathbb{R}^N)}\ge \frac{1}{T_{N,1/2}}\,\|u\|^N_{L^{2\,N}(\Omega)},
\]
where we used the definition \eqref{talentiana}.
Then by joining the two previous estimates, appealing to the definition of $\alpha^s_N(\Omega)$ and recalling that $s\ge 3/4$, we get
\[
\alpha^s_N(\Omega)\ge \frac{1}{C}\,\dfrac{\Big(\mathrm{diam}(\Omega)\Big)^{N\,\frac{1-2\,s}{2}}}{s\,(1-s)},
\]
where $C=C(N)>0$.
By inserting this estimate in \eqref{alpelo}, we get the conclusion in this case as well.
\vskip.2cm\noindent
For $s\,p>N$, we already know that $\widetilde W^{s,p}_0(\Omega)\hookrightarrow C^{0,s-N/p}$, but of course we need to estimate the embedding constant in terms of $s$. We take $x_0\in \mathbb{R}^N$ and $R>0$. We consider the ball $B_R(x_0)$ having radius $R$ and centered at $x_0$, then we have
\[
\begin{split}
\int_{B_R(x_0)}\left|u-\fint_{B_R(x_0)} u\,dz\right|^p\,dx&\le \int_{B_R(x_0)}\,\fint_{B_{R}(x_0)} |u(x)-u(z)|^p\,dz\,dx\\
\end{split}
\]
We now observe that
\[
\begin{split}
\int_{B_R(x_0)}\fint_{B_{R}(x_0)} |u(x)-u(z)|^p\,dz\,dx&\le \int_{B_R(x_0)}\left(\fint_{B_{2\,R}(x)} |u(x)-u(z)|^p\,dz\right)dx\\
&=\int_{B_R(x_0)}\left(\fint_{B_{2\,R}(0)} |u(x)-u(x+h)|^p\,dh\right)dx,
\end{split}
\]
so that by exchanging the order of integration in the last integral
\[
\int_{B_R(x_0)}\left|u-\fint_{B_R(x_0)} u\,dz\right|^p\,dx\le C\,(2\,R)^{s\,p}\sup_{0<|h|<2\,R} \int_{B_R(x_0)} \frac{|u(x)-u(x+h)|^p}{|h|^{s\,p}}\,dx,
\]
for $C=C(N)>0$.
If we now divide by $R^{N}$ and use once again \cite[Lemma A.1]{BraLinPar}, we get
\[
\left(\fint_{B_R(x_0)}\left|u-\fint_{B_R(x_0)} u\,dz\right|^p\,dx\right)^\frac{1}{p}\le C\,(2\,R)^{s-\frac{N}{p}}\,(1-s)^\frac{1}{p}\,[u]_{W^{s,p}(\mathbb{R}^N)}.
\]
By arbitrariness of $R$ and $x_0$, we obtain that $u$ is in $C^{0,s-N/p}(\mathbb{R}^N)$ by Campanato's Theorem (see \cite[Theorem 2.9]{Gi}), with the estimate
\begin{equation}
\label{campanella}
|u(x)-u(y)|\le C\, (1-s)^\frac{1}{p}\,[u]_{W^{s,p}(\mathbb{R}^N)}\,|x-y|^{s-\frac{N}{p}},
\end{equation}
where $C=C(N,p)>0$. The last estimate is true for every $u\in \widetilde W^{s,p}_0(\Omega)$. On the other hand, if $u\in \widetilde W^{s,p}_0(\Omega)$ is an eigenfunction with eigenvalue $\lambda$, then by the equation we also have
\[
[u]_{W^{s,p}(\mathbb{R}^N)}=\lambda^\frac{1}{p}\,\|u\|_{L^p(\Omega)}.
\]
By inserting this estimate in \eqref{campanella}, we get
\begin{equation}
\label{holder}
|u(x)-u(y)|\le \Big[C\,(1-s)\,\lambda\Big]^\frac{1}{p}\,\|u\|_{L^p(\Omega)}\,|x-y|^{s-\frac{N}{p}},\qquad x,y\in\mathbb{R}^N.
\end{equation}
Finally, the estimate \eqref{bounder} follows from \eqref{holder} by taking $y\in \mathbb{R}^N\setminus\Omega$.
\end{proof}
\begin{remark}
Though we will not need it here, for the conformal case $s\,p=N$ a global $L^\infty$ estimate can be found in \cite[Theorem 3.3]{BraLinPar}.
\end{remark}

\section{A $\Gamma-$convergence result}
\label{sectGamma}

In this section we will prove the following result.
\begin{theorem}[$\Gamma-$convergence]
\label{thm:gamma}
Let $1<p<\infty$ and $\Omega\subset\mathbb{R}^N$ be an open and bounded set, with Lipschitz boundary. We consider $\{s_k\}_{h\in\mathbb{N}}$ a sequence of strictly increasing positive numbers, such that
$s_k$ goes to $1$ as $k$ goes to $\infty$. 
Then
\begin{equation}
\label{gamma0}
\mathcal E_{1,p}(u)=\Big(\Gamma-\lim_{k\to\infty}\mathcal E_{s_k,p}\Big)(u),
\qquad\mbox{ for all } u\in L^p(\Omega).
\end{equation}
where $\mathcal{E}_{s_k,p}$ and $\mathcal{E}_{1,p}$ are the functionals defined by \eqref{funzionales} and \eqref{funzionale1}.
\end{theorem}
This $\Gamma-$convergence result will follow from Propositions \ref{prop:limsup} and \ref{prop:liminf} below. Before proceeding further with the proof of this result, let us highlight that by combining Theorem \ref{thm:gamma} and \cite[Proposition 6.25]{dalmaso1993}, we get the following.
\begin{corollary}
\label{cor:gammacor}
Under the assumptions of Theorem \ref{thm:gamma}, we also consider a sequence of functions $\{F_{s_k}\}_{k\in\mathbb{N}}\subset L^{p'}(\Omega)$ weakly converging in $L^{p'}(\Omega)$ to $F$. If we introduce the functionals defined on $L^p(\Omega)$ by
\begin{equation}
\label{Fps}
\mathcal F_{s_k,p}(u):=\mathcal{E}_{s_k,p}(u)^p+\int_\Omega F_{s_k}\, u\,dx\qquad \mbox{ and }\qquad \mathcal F_{1,p}(u):=\mathcal{E}_{1,p}(u)^p+\int_\Omega F\,u\,dx,
\end{equation}  
then we also have
\[
\mathcal F_{1,p}(u)=\Big(\Gamma-\lim_{h\to\infty}\mathcal F_{s_k,p}\Big)(u),
\qquad\mbox{ for all } u\in L^p(\Omega).
\]
\end{corollary}

\subsection{The $\Gamma-\limsup$ inequality}

\begin{proposition}[$\Gamma-\limsup$ inequality]
\label{prop:limsup}
Let $u\in L^p(\Omega)$ and let $\{s_k\}_{k\in\mathbb{N}}$ be a sequence of strictly increasing positive numbers, such that $s_k$ converges to $1$ as $k$ goes to $\infty$. Then there exists a sequence $\{u_k\}_{k\in\mathbb{N}}\subset \widetilde W^{s_k,p}_0(\Omega)$ such that
\[
\limsup_{k\to\infty} \mathcal{E}_{s_k,p}(u_k)^p\le \mathcal{E}_{1,p}(u)^p.
\]
\end{proposition}
\begin{proof}
If $u\not\in W^{1,p}_0(\Omega)$, there is nothing to prove, thus let us take $u\in W^{1,p}_0(\Omega)$.
If we take the constant sequence $u_k=u$ and then apply the modification of Bourgain-Brezis-Mironescu result of Proposition \ref{norm-cont}, we obtain
\[
\limsup_{h\to\infty} (1-s_k)\,[u_k]^p_{W^{s_k,p}(\mathbb{R}^N)}=K(p,N)\,\int_\Omega |\nabla u|^p\,dx=\mathcal{E}_{1,p}(u)^p,
\]
concluding the proof.
\end{proof}
In order to prove the $\Gamma-\liminf$ inequality, we need to find a different characterization of the constant $K(p,N)$. The rest of this subsection is devoted to this issue.
\par
In what follows, we note by $Q=(-1/2,1/2)^N$ the open $N-$dimensional cube of side length $1$. Given $a \in \mathbb{R}^N$, we define the linear function $\Psi_a(x)=\langle a,x\rangle$. For every $a\in\mathbb{S}^{N-1}$, we define the constant
\begin{equation}
\label{gammapna}
\Theta(p,N;a) := \inf \left\{ \liminf_{s \nearrow 1} (1-s)\,[u_s]_{W^{s,p}(Q)}^p\,:\, u_s \to \Psi_a \text{ in }L^p(Q)\right\}.
\end{equation}
We will show in Lemma \ref{lm:costanti} that indeed this quantity does not depend on the direction $a$.
\begin{remark}
If $a\in \mathbb{R}^N$ with $|a|\not =0$, then we have
\begin{equation}
\label{nonunit}
\inf \left\{ \liminf_{s \nearrow 1} (1-s)\,[u_s]_{W^{s,p}(Q)}^p\,:\, u_s \to \Psi_a \text{ in }L^p(Q)\right\}=|a|^p\, \Theta\left(p,N;\frac{a}{|a|}\right).
\end{equation}
\end{remark}
\begin{remark}
For every $1<p<\infty$ and every $a\in\mathbb{S}^{N-1}$ we have
\begin{equation}
\label{unverso}
\Theta(p,N;a)\le K(p,N),
\end{equation}
where $K(p,N)$ is the constant defined in \eqref{defK}. Indeed, by definition of $\Theta(p,N;a)$, if we take the constant sequence $u_s=\Psi_a$ and use the Bourgain-Brezis-Mironescu result, we get
\[
\Theta(p,N;a)\le \liminf_{s \nearrow 1} (1-s)\,[\Psi_a]_{W^{s,p}(Q)}^p=K(p,N)\,\int_Q |\nabla \Psi_a|^p\,dx.
\]
This proves \eqref{unverso}, since $\nabla \Psi_a$ has unit norm in $L^p(Q)$.
\end{remark}
We are going to prove that indeed $K(p,N)=\Theta(p,N;a)$ for every $a\in\mathbb{S}^{N-1}$. To this aim, we first need a couple of technical results. In what follows, by $W^{s,p}_0(Q)$ we note the completion of $C^\infty_0(Q)$ with respect to the semi-norm 
\[
[u]_{W^{s,p}(Q)}:=\left(\int_Q\int_Q \frac{|u(x)-u(y)|^p}{|x-y|^{N+s\,p}}\,dx\,dy\right)^\frac{1}{p}.
\]
\begin{lemma}
For every $1<p<\infty$ and every $a\in\mathbb{S}^{N-1}$ we have
\begin{equation}
\label{bordouguali}
\Theta(p,N;a) = \inf \left\{ \liminf_{s\nearrow 1} (1-s)\,[v_s]_{W^{s,p}(Q)}^p: \begin{array}{c}v_s \to \Psi_a \text{ in }L^p(Q)\\v_s-\Psi_a\in W^{s,p}_0(Q) \end{array}\right\}.
\end{equation}
\end{lemma}
\begin{proof}
By the definition \eqref{gammapna} of $\Theta(p,N;a)$, we already know that
\[
\Theta(p,N;a) \le \inf \left\{ \liminf_{s\nearrow 1} (1-s)\,[v_s]_{W^{s,p}(Q)}^p\,: \begin{array}{c}v_s \to \Psi_a \text{ in }L^p(Q)\\v_s-\Psi_a\in W^{s,p}_0(Q) \end{array} \right\}.
\]
In order to prove the reverse inequality, let us take a sequence $\{s_k\}_{k\in\mathbb{N}}$ such that $0<s_k<1$ and $s_k\nearrow 1$. Then we take $\{v_{k}\}_{k\in\mathbb{N}}$ such that
\[
(1-s_k)\,[v_k]^p_{W^{s_k,p}(Q)}<+\infty \qquad \mbox{ and }\qquad \lim_{k\to\infty}\|v_k-\Psi_a\|_{L^p(\Omega)}=0.
\]
Without loss of generality, we can assume that $s_k\,p>1$, so that for the space $W^{s_k,p}_0(\Omega)$ we have the Poincar\'e inequality \eqref{spoinbis}.
We introduce a smooth cut-off function $\eta\in C^\infty_0(Q)$ such that
\[
0\le \eta\le 1,\qquad \eta\equiv 1,\mbox{ on } \tau\,Q,\qquad |\nabla \eta|\le \frac{10}{1-\tau},
\]
for some parameter $0<\tau<1$. Then we define the sequence $\{w_k\}_{k\in\mathbb{N}}$ by
\[
w_k:=v_k\,\eta+\Psi_a\,(1-\eta).
\]
We observe that by construction we have $w_k-\Psi_a\in W^{s,p}_0(Q)$. Moreover we have
\[
\begin{split}
\int_Q |w_k-\Psi_a|^p\,dx=\int_Q \eta^p\,|v_k-\Psi_a|^p\,dx&\le \int_Q |v_k-\Psi_a|^p\,dx,
\end{split}
\]
thus $w_k$ still converges in $L^p(Q)$ to $\Psi_a$. We now have to estimate the Gagliardo semi-norm of $w_k$. To this aim, we first observe that
\begin{equation}
\label{ciao}
\begin{split}
w_k(x)-w_k(y)&=v_k(x)\,\eta(x)+\Psi_a(x)\,(1-\eta(x))-v_k(y)\,\eta(y)-\Psi_a(y)\,(1-\eta(y))\\
&=\eta(x) \Big(v_k(x)-v_k(y)\Big)+(1-\eta(x))\Big(\Psi_a(x)-\Psi_a(y)\Big)\\
&+\Big(\eta(x)-\eta(y)\Big)\,\Big(v_k(y)-\Psi_a(y)\Big).
\end{split}
\end{equation}
Let us set
\[
V(x,y)=\eta(x)\,\frac{v_k(x)-v_k(y)}{|x-y|^{\frac{N}{p}+s_k}}\qquad \Pi(x,y)=(1-\eta(x))\,\frac{\Psi_a(x)-\Psi_a(y)}{|x-y|^{\frac{N}{p}+s_k}},
\]
and
\[
Z(x,y)=\frac{\eta(x)-\eta(y)}{|x-y|^{\frac{N}{p}+s_k}}\,\Big(v_k(y)-\Psi_a(y)\Big).
\]
Then by definition of $w_k$, \eqref{ciao} and Minkowski inequality we have
\[
\begin{split}
[w_k]_{W^{s_k,p}(Q)}&=\|V+\Pi+Z\|_{L^p(Q\times Q)}\le \|V\|_{L^p(Q\times Q)}+\|\Pi\|_{L^p(Q\times Q)}+\|Z\|_{L^p(Q\times Q)}\\
&=\left(\int_Q \int_Q \frac{|v_k(x)-v_k(y)|^p}{|x-y|^{N+s_k\,p}}\,\eta(x)^p\,dx\,dy\right)^\frac{1}{p}\\
&+\left(\int_Q \int_Q \frac{|\Psi_a(x)-\Psi_a(y)|^p}{|x-y|^{N+s_k\,p}}\,(1-\eta(x))^p\,dx\,dy\right)^\frac{1}{p}\\
&+\left(\int_Q \int_Q \frac{|\eta(x)-\eta(y)|^p}{|x-y|^{N+s_k\,p}}\,|v_k(y)-\Psi_a(y)|^p\,dx\,dy\right)^\frac{1}{p}.
\end{split}
\]
By using the properties of $\eta$, we have obtained
\begin{equation}
\label{better}
\begin{split}
[w_k]_{W^{s_k,p}(Q)}&\le [v_k]_{W^{s_k,p}(Q)}+\left(\int_Q \int_{Q\setminus\tau\,Q} \frac{|\Psi_a(x)-\Psi_a(y)|^p}{|x-y|^{N+s_k\,p}}\,dx\,dy\right)^\frac{1}{p}\\
&+\frac{C}{1-\tau}\,\left(\int_Q \left(\int_Q \frac{dx}{|x-y|^{N+s_k\,p-p}}\right)\,|v_k(y)-\Psi_a(y)|^p\,dy\right)^\frac{1}{p}.
\end{split}
\end{equation}
We have to estimate the last two integrals. By recalling that $\Psi_a(x)=\langle a,x\rangle$, we have\footnote{We use that for $x\in Q$
\[
\int_Q \frac{1}{|x-y|^{N+s_k\,p-p}}\,dy\le\int_{B_{\sqrt{N}}(x)} \frac{1}{|x-y|^{N+s_k\,p-p}}\,dy=\int_{B_{\sqrt{N}}(0)} \frac{1}{|y|^{N+s_k\,p-p}}\,dy=\frac{C}{1-s_k},
\]
with $C=C(N,p)>0$}
\[
\begin{split}
\int_Q \int_{Q\setminus\tau\,Q} \frac{|\Psi_a(x)-\Psi_a(y)|^p}{|x-y|^{N+s_k\,p}}\,dx\,dy&\le \int_{Q\setminus\tau\,Q} \int_Q \frac{1}{|x-y|^{N+s_k\,p-p}}\,dy\,dx\\
&\le \frac{C}{1-s_k}\,|Q\setminus\tau\,Q|.
\end{split}
\]
For the other integral, we have
\[
\left(\int_Q \left(\int_Q \frac{dx}{|x-y|^{N+s_k\,p-p}}\right)\,|v_k(y)-\Psi_a(y)|^p\,dy\right)^\frac{1}{p}\le \frac{C}{(1-s_k)^\frac{1}{p}}\,\|v_k-\Psi_a\|_{L^p(\Omega)},
\]
with $C=C(N,p)>0$. 
By collecting all these estimates and using them in \eqref{better}, we get
\[
\begin{split}
\liminf_{k\to\infty}(1-s_k)^\frac{1}{p}\,[w_k]_{W^{s_k,p}(Q)}&\le \liminf_{k\to\infty}(1-s_k)^\frac{1}{p}\,[v_k]_{W^{s_k,p}(Q)}\\
&+\frac{C}{1-\tau}\,\liminf_{k\to\infty} \|v_k-\Psi_a\|_{L^p(Q)}+C\,|Q\setminus \tau\, Q|\\
&=\liminf_{k\to\infty}(1-s_k)^\frac{1}{p}\,[v_k]_{W^{s_k,p}(Q)}+C\,|Q\setminus \tau\, Q|.\\
\end{split}
\]
By arbitrariness of $0<\tau<1$, this finally proves the desired result.
\end{proof}
Before proceeding further, we need to know that linear functions are $(s,p)-$harmonic. 
\begin{lemma}[Linear functions are $(s,p)-$harmonic]
\label{lm:lineari}
Let $a \in \mathbb{S}^{N-1}$, we define $\Psi_a(x)=\langle a,x\rangle$ as above. Let $1<p<\infty$ and $s\in(0,1)$ such that $s>(p-1)/p$. Then for every $\varphi\in C^\infty_0(Q)$ we have
\begin{equation}
\label{lineari}
\int_{\mathbb{R}^N}\int_{\mathbb{R}^N} \frac{|\Psi_a(x)-\Psi_a(y)|^{p-2}\,(\Psi_a(x) - \Psi_a(y)) }{|x-y|^{N+s\,p}}\,(\varphi(x)-\varphi(y))\,dx\,dy=0.
\end{equation}
\end{lemma}
\begin{proof}
We first observe that the double integral is well-defined and absolutely convergent.
Indeed,
\[
\begin{split}
\int_{\mathbb{R}^N}\int_{\mathbb{R}^N} &\frac{|\Psi_a(x)-\Psi_a(y)|^{p-2}\,(\Psi_a(x) - \Psi_a(y)) }{|x-y|^{N+s\,p}}\,(\varphi(x)-\varphi(y))\,dx\,dy\\
&=\int_{Q}\int_{Q} \frac{|\Psi_a(x)-\Psi_a(y)|^{p-2}\,(\Psi_a(x) - \Psi_a(y)) }{|x-y|^{N+s\,p}}\,(\varphi(x)-\varphi(y))\,dx\,dy\\
&+2\, \int_{Q}\int_{\mathbb{R}^N\setminus Q}\frac{|\Psi_a(x)-\Psi_a(y)|^{p-2}\,(\Psi_a(x) - \Psi_a(y)) }{|x-y|^{N+s\,p}}\,\varphi(x)\,dx\,dy.
\end{split}
\]
For the first term we have
\[
\int_{Q}\int_{Q} \frac{|\Psi_a(x)-\Psi_a(y)|^{p-1}}{|x-y|^{N+s\,p}}\,|\varphi(x)-\varphi(y)|\,dx\,dy\le \int_{Q}\int_Q \frac{\|\nabla \varphi\|_{L^\infty} \,dx\,dy}{|x-y|^{N+s\,p-p}}<+\infty.
\]
For the second one, by observing that the integral in the $x$ variable is equivantely performed on $K:=\mathrm{spt}(\varphi)\Subset Q$, we get
\[
\int_{Q}\int_{\mathbb{R}^N\setminus Q}\frac{|\Psi_a(x)-\Psi_a(y)|^{p-1}}{|x-y|^{N+s\,p}}\,|\varphi(x)|\,dx\,dy\le \int_{K}\int_{\mathbb{R}^N\setminus Q}\frac{\|\varphi\|_{L^\infty(Q)}\,dx\,dy}{|x-y|^{N+s\,p-p+1}}<+\infty,
\]
provided that $s>(p-1)/p$.
\par
In order to prove \eqref{lineari}, for every $\varepsilon>0$ we have
\[
\begin{split} 
\int_{\mathbb{R}^N}\int_{\mathbb{R}^N} &\frac{|\Psi_a(x)-\Psi_a(y)|^{p-2}\,(\Psi_a(x) - \Psi_a(y)) }{|x-y|^{N+s\,p}}\,(\varphi(x)-\varphi(y))\,dx\,dy \\
&=\iint_{\{|x-y|\ge \varepsilon\}} \frac{|\Psi_a(x)-\Psi_a(y)|^{p-2}\,(\Psi_a(x) - \Psi_a(y)) }{|x-y|^{N+s\,p}}\,\varphi(x)\,dx\,dy\\
&-\iint_{\{|x-y|\ge \varepsilon\}} \frac{|\Psi_a(x)-\Psi_a(y)|^{p-2}\,(\Psi_a(x) - \Psi_a(y)) }{|x-y|^{N+s\,p}}\,\varphi(y)\,dx\,dy\\
&+\iint_{\{|x-y|<\varepsilon\}} \frac{|\Psi_a(x)-\Psi_a(y)|^{p-2}\,(\Psi_a(x) - \Psi_a(y)) }{|x-y|^{N+s\,p}}\,(\varphi(x)-\varphi(y))\,dx\,dy\\
&=2\,\iint_{\{|x-y|\ge \varepsilon\}} \frac{|\Psi_a(x)-\Psi_a(y)|^{p-2}\,(\Psi_a(x) - \Psi_a(y)) }{|x-y|^{N+s\,p}}\,\varphi(x)\,dx\,dy\\
&+\iint_{\{|x-y|<\varepsilon\}} \frac{|\Psi_a(x)-\Psi_a(y)|^{p-2}\,(\Psi_a(x) - \Psi_a(y)) }{|x-y|^{N+s\,p}}\,(\varphi(x)-\varphi(y))\,dx\,dy\\
&=2\,\int_{\mathbb{R}^N} \left(\int_{\{|h|>\varepsilon\}}\frac{|\langle a,h\rangle|^{p-2}\,\langle a,h\rangle}{|h|^{N+s\,p}}\,dh\right)\,\varphi(x)\,dx\\
&+\iint_{\{|x-y|<\varepsilon\}} \frac{|\Psi_a(x)-\Psi_a(y)|^{p-2}\,(\Psi_a(x) - \Psi_a(y)) }{|x-y|^{N+s\,p}}\,(\varphi(x)-\varphi(y))\,dx\,dy\\
&=\iint_{\{|x-y|<\varepsilon\}} \frac{|\Psi_a(x)-\Psi_a(y)|^{p-2}\,(\Psi_a(x) - \Psi_a(y)) }{|x-y|^{N+s\,p}}\,(\varphi(x)-\varphi(y))\,dx\,dy,
\end{split}
\]
where we used that by symmetry
\[
\int_{\{|h|>\varepsilon\}}\frac{|\langle a,h\rangle|^{p-2}\,\langle a,h\rangle}{|h|^{N+s\,p}}\,dh=0.
\]
Moreover, we have (we still denote $K=\mathrm{spt}(\varphi)$)
\[
\begin{split}
\Big|\iint_{\{|x-y|<\varepsilon\}}& \frac{|\Psi_a(x)-\Psi_a(y)|^{p-2}\,(\Psi_a(x) - \Psi_a(y)) }{|x-y|^{N+s\,p}}\,(\varphi(x)-\varphi(y))\,dx\,dy\Big|\\
&\le\|\nabla \varphi\|_{L^\infty}\, \int_{K+B_\varepsilon(0)}\left(\int_{\{|h|<\varepsilon\}} \frac{1}{|h|^{N+s\,p-p}}\,dh\right)\,dx\\
&\le \frac{C\, \varepsilon^{p\,(1-s)}}{1-s}\,\|\nabla \varphi\|_{L^\infty}\,\Big|K+B_1(0)\Big|.
\end{split}
\]
By arbitrariness of $\varepsilon$ we get the conclusion. 
\end{proof}

\begin{lemma}
\label{lm:minimizzatore}
Let $a \in \mathbb{S}^{N-1}$. For every $1<p<\infty$ and $s\in(0,1)$ such that $s\,p>1$, let $u_s$ be the unique solution of
\begin{equation} 
\label{minimalityus} 
\min\Big\{[v]^p_{W^{s,p}(Q)} \,:\, v-\Psi_a\in W^{s,p}_0(Q) \Big\}.
\end{equation}
Then, $u_s$ converges to $\Psi_a$ in $L^p(\Omega)$ as $s$ goes to $1$. Moreover, we have
\begin{equation}
\label{limites}
K(p,N)=\lim_{s \nearrow 1} (1-s)\,[\Psi_a]_{W^{s,p}(Q)}^p =\lim_{s \nearrow 1} (1-s)\,[u_s]_{W^{s,p}(Q)}^p
\end{equation}
\end{lemma}
\begin{proof}
Since we are interested in the limit as $s$ goes to $1$, without loss of generality we can further assume that $s>(p-1)/p$ as well, i.e.
\[
s>\max\left\{\frac{1}{p},\frac{p-1}{p}\right\}.
\]
The existence of a (unique, by strict convexity) solution $u_s$ follows by the Direct Methods, since coercivity of the functional $v\mapsto [v]^p_{W^{s,p}(Q)}$ can be inferred thanks to Poincar\'e inequality \eqref{spoinbis} (here we use the assumption $s\,p>1$).
We take $\varphi \in W^{s,p}_0(Q)$, by minimality of $u_s$ there holds
\begin{equation}
\label{eguazione}
\int_Q \int_Q \frac{|u_s(x)-u_s(y)|^{p-2}(u_s(x) - u_s(y))}{|x-y|^{N+s\,p}}\,(\varphi(x)-\varphi(y))\,dx\,dy = 0.
\end{equation}
Still by minimality of $u_s$, we also get
\begin{equation}
\label{armonia}
[u_s]_{W^{s,p}(Q)}\le [\Psi_a]_{W^{s,p}(Q)},
\end{equation}
since $\Psi_a$ is admissible for problem \eqref{minimalityus}.
On the other hand, the linear function $\Psi_a$ is ``almost'' a solution of \eqref{minimalityus}, thanks to Lemma \ref{lm:lineari}. Indeed from \eqref{lineari}, for every $\varphi\in C^\infty_0(Q)$ we get
\[
\begin{split}
\int_{Q} \int_{Q} &\frac{|\Psi_a(x)-\Psi_a(y)|^{p-2}(\Psi_a(x) - \Psi_a(y)) }{|x-y|^{N+s\,p}}(\varphi(x)-\varphi(y))\,dx\,dy \\ & = -  \int_{Q} \int_{\mathbb{R}^N\setminus Q} \frac{|\Psi_a(x)-\Psi_a(y)|^{p-2}(\Psi_a(x) - \Psi_a(y)) }{|x-y|^{N+s\,p}}(\varphi(x)-\varphi(y))\,dx\,dy \\ & - \int_{\mathbb{R}^N\setminus Q} \int_{Q} \frac{|\Psi_a(x)-\Psi_a(y)|^{p-2}(\Psi_a(x) - \Psi_a(y)) }{|x-y|^{N+s\,p}}(\varphi(x)-\varphi(y))\,dx\,dy \\ 
& = -2\,\int_{\mathrm{spt}(\varphi)} \int_{\mathbb{R}^N\setminus Q} \frac{|\Psi_a(x)-\Psi_a(y)|^{p-2}\,(\Psi_a(x) - \Psi_a(y)) }{|x-y|^{N+s\,p}}\,\varphi(x)\,dx\,dy.
\end{split}
\]
Thus we obtain
\[
\begin{split}
 \bigg| \int_{Q} \int_{Q} &\frac{|\Psi_a(x)-\Psi_a(y)|^{p-2}\,(\Psi_a(x) - \Psi_a(y)) }{|x-y|^{N+s\,p}}\,(\varphi(x)-\varphi(y))\,dx\,dy \bigg| \\ 
 & \leq 2\,\int_{Q} \int_{\mathbb{R}^N\setminus Q} \frac{|\varphi(x)|}{|x-y|^{N+s\,p-p+1}}\,dx\,dy\\
 &\le 2\, \int_{Q} \int_{\mathbb{R}^N\setminus B_{\delta_Q(x)}(x)} \frac{|\varphi(x)|}{|x-y|^{N+s\,p-p+1}}\,dx\,dy,
\end{split}
\]
where as before we set $\delta_Q(x)=\text{dist}(x,\partial Q)$. Hence,
\[
\begin{split}
\bigg| \int_{Q} \int_{Q} &\frac{|\Psi_a(x)-\Psi_a(y)|^{p-2}(\Psi_a(x) - \Psi_a(y)) }{|x-y|^{N+s\,p}}(\varphi(x)-\varphi(y))\,dx\,dy \bigg|\\
&  \leq 2\,N\omega_N\, \int_{Q} \left(\int_{\delta_Q(x)}^{+\infty} \frac{1}{\varrho^{2+s\,p-p}}\,d\varrho\right)\,|\varphi(x)|\,dx \\ & = \frac{2\,N\omega_N}{1+s\,p-p} \int_Q \frac{|\varphi|}{\delta_Q^{1+s\,p-p}}\,dx \\ 
& \leq \frac{2\,N\omega_N}{1+s\,p-p} \left(\int_Q \frac{|\varphi|^p}{\delta_Q^{s\,p}}\,dx\right)^{\frac{1}{p}} \left( \int_Q \delta^{p\,(1-s)}_Q\,dx \right)^{\frac{p-1}{p}} \leq C \left(\int_Q \frac{|\varphi|^p}{\delta^{s\,p}_Q}\,dx\right)^{\frac{1}{p}}.
\end{split}
\]
Since we are assuming $s > 1/p$, we can apply Hardy inequality \eqref{hardysubopt}
to the last term and obtain 
\begin{equation} 
\label{estimatelinear} 
\begin{split}
\int_{Q} \int_{Q} \frac{|\Psi_a(x)-\Psi_a(y)|^{p-2}\,(\Psi_a(x) - \Psi_a(y)) }{|x-y|^{N+s\,p}}&(\varphi(x)-\varphi(y))\,dx\,dy\\
&\leq \frac{C}{s\,p-1}\,[\varphi]_{W^{s,p}(Q)},
\end{split}
\end{equation}
for some constant $C=C(N,p)>0$ (observe that we used that $1-s<1$).
From \eqref{eguazione} and \eqref{estimatelinear} we finally obtain for every $\varphi\in C^\infty_0(Q)$
\begin{equation} 
\label{stimaconvarphi} 
\begin{split}
\int_{Q} \int_{Q} &\frac{|\Psi_a(x)-\Psi_a(y)|^{p-2}\,(\Psi_a(x) - \Psi_a(y)) - |u_s(x)-u_s(y)|^{p-2}(u_s(x) - u_s (y)) }{|x-y|^{N+sp}}\\
&\times(\varphi(x)-\varphi(y))\,dx\,dy\leq C\,[\varphi]_{W^{s,p}(Q)}
\end{split}
\end{equation}
By density, the previous estimate is still true for every $\varphi \in W^{s,p}_0(Q)$, thus we can use \eqref{stimaconvarphi} with $\varphi = \Psi_a - u_s$. We distinguish two cases. 
\vskip.2cm\noindent
\underline{\it Case $p \geq 2$}. We use the basic inequality $(|s|^{p-2}s - |t|^{p-2}t)(s-t) \geq 2^{2-p}|s-t|^p$ in order to obtain from \eqref{stimaconvarphi}
\[ 
[\Psi_a - u_s]_{W^{s,p}(Q)}^p \leq C\,[\Psi_a - u_s]_{W^{s,p}(Q)},
\]
for a constant $C=C(N,p)>0$.
This implies
\begin{equation}
\label{2p}
\lim_{s\nearrow 1} (1-s)\,[\Psi_a - u_s]^p_{W^{s,p}(Q)}=0. 
\end{equation}
\underline{\it Case $1<p<2$}. We use the inequality 
\[
(|s|^{p-2}s - |t|^{p-2}t)\,(s-t) \geq (p-1)\frac{|s-t|^2}{(|s|^2+|t|^2)^\frac{2-p}{2}}
\]
which gives
\[
|s-t|^p \leq \left(\frac{1}{p-1}\right)^{\frac{p}{2}} \Big[(|s|^{p-2}s - |t|^{p-2}t)(s-t)\Big]^{\frac{p}{2}} (|s|^2+|t|^2)^{\frac{2-p}{2}\frac{p}{2}} .
\]
We set for notational simplicity $U_s(x,y)=u_s(x)-u_s(y)$ and $U(x,y)=\Psi_a(x)-\Psi_a(y)$. Then,
\[
\begin{split} 
[\Psi_a&-u_s]^p_{W^{s,p}(Q)}\\
 & \le \left(\frac{1}{p-1}\right)^{\frac{p}{2}} \int_Q \int_Q \frac{\Big[(|U|^{p-2}\,U - |U_s|^{p-2}\,U_s)(U-U_s)\Big]^{\frac{p}{2}} (|U|^2+|U_s|^2)^{\frac{2-p}{2}\frac{p}{2}}}{|x-y|^{N+sp}} \\ & \leq \left(\frac{1}{p-1}\right)^{\frac{p}{2}} \left( \int_Q \int_Q \frac{(|U|^{p-2}\,U - |U_s|^{p-2}\,U_s)(U-U_s)}{|x-y|^{N+sp}}\right)^{\frac{p}{2}}\\
&\times \left( \int_Q \int_Q \frac{(|U|^2+|U_s|^2)^{\frac{p}{2}}}{|x-y|^{N+sp}}\right)^{\frac{2-p}{2}} \\ & \leq C\, \Big([\Psi_a-u_s]_{W^{s,p}(Q)}\Big)^\frac{p}{2}\, \left( [u_s]^p_{W^{s,p}(Q)}+[\Psi_a]_{W^{s,p}(Q)}^p\right)^{\frac{2-p}{2}}, 
\end{split}
\]
where we used H\"{o}lder inequality with exponents $2/p$ and $2/(2-p)$, relation \eqref{stimaconvarphi} and the subadditivity of the function $t \mapsto t^{p/2}$. The previous estimate and \eqref{armonia} imply
\[ 
[\Psi_a-u_s]^p_{W^{s,p}(Q)} \leq C\,[\Psi_a]_{W^{s,p}(Q)}^{(2-p)\,p}, 
\]
that is
\[
(1-s)\,[\Psi_a-u_s]^p_{W^{s,p}(Q)} \leq C\,(1-s)^{p-1}\,\Big((1-s)\,[\Psi_a]^p_{W^{s,p}(Q)}\Big)^{{2-p}}.
\]
Since $2-p<1$, after a simplification we get
\[
\left((1-s)\,[\Psi_a-u_s]^p_{W^{s,p}(Q)}\right)^{p-1} \leq C\,(1-s)^{p-1}.
\]
It thus follows again
\begin{equation}
\label{p2} 
\lim_{s\nearrow 1} (1-s)\, [\Psi_a - u_s]^p_{W^{s,p}(Q)}=0.
\end{equation}
Observe that as a byproduct of \eqref{2p} and \eqref{p2}, we also get
\[
\begin{split}
\lim_{s\nearrow 1} \Big|(1-s)^{\frac{1}{p}}\,[\Psi_a]_{W^{s,p}(Q)}&-(1-s)^{\frac{1}{p}}\,[u_s]_{W^{s,p}(Q)}\Big|\\
&\le\lim_{s\nearrow 1} (1-s)^\frac{1}{p}\, [\Psi_a - u_s]_{W^{s,p}(Q)}=0.
\end{split}
\]
This shows that 
\[
\lim_{s\nearrow 1} (1-s)^\frac{1}{p}\, [\Psi_a]_{W^{s,p}(Q)}=\lim_{s\nearrow 1} (1-s)^\frac{1}{p}\, [u_s]_{W^{s,p}(Q)},
\]
thus \eqref{limites} is proved.
\vskip.2cm\noindent
Finally, since $u_s-u\in W^{s,p}_0(Q)$, $Q$ is a convex set and we are assuming $s\,p>1$, we can use Poincar\'e inequality \eqref{spoinbis} in conjuction with \eqref{2p} or \eqref{p2}. In both cases we have
\[ 
\lim_{s\nearrow 1}\|u_s - \Psi_a\|^p_{L^p(Q)} \leq \lim_{s\nearrow 1} \frac{C}{s\,p-1}\,(1-s)\, [u_s - \Psi_a]^p_{W^{s,p}(Q)}=0,
\]
where $C=C(N,p)>0$.
This concludes the proof.
\end{proof}
Finally, we can prove an equivalent characterization of $K(p,N)$.
\begin{lemma}
\label{lm:costanti}
Let $1<p<\infty$ and $a\in\mathbb{S}^{N-1}$, then we have
\[
K(p,N)=\Theta(p,N;a).
\]
In particular, $\Theta(p,N;a)$ does not depend on the direction $a$.
\end{lemma}
\begin{proof}
By \eqref{unverso} we know that $K(p,N)\ge \Theta(p,N;a)$. In order to prove the reverse inequality, we define the linear function $\Psi_a(x)=\langle a,x\rangle$. 
Let $v_s\in W^{s,p}(Q)$ be a sequence converging to $\Psi_a$ in $L^p(Q)$ and such that $v_s-\Psi_a\in W^{s,p}_0(Q)$. We consider the function $u_s$ defined in Lemma \ref{lm:minimizzatore}, then from \eqref{limites} we get
\[
\liminf_{s\nearrow 1}(1-s)\,[v_s]^p_{W^{s,p}(Q)}\ge \liminf_{s\nearrow 1} (1-s)\,[u_s]^p_{W^{s,p}(Q)}=K(p,N).
\]
By appealing to \eqref{bordouguali}, we get $\Theta(p,N;a)=K(p,N)$. 
\end{proof}

\subsection{The $\Gamma-\liminf$ inequality}

At first, we need a technical result which will be used various times.
\begin{lemma}
\label{lm:utile}
Let $1<p<\infty$ and $s_0\in(0,1)$. Let $\Omega\subset\mathbb{R}^N$ be an open and bounded set with Lipschitz boundary. For every family of functions $\{u_s\}_{s\in(s_0,1)}$ such that $u_s\in \widetilde W^{s,p}_0(\Omega)$ and
\begin{equation}
\label{ipotesi_utile}
(1-s)\,[u_s]^p_{W^{s,p}(\mathbb{R}^N)}\le L,
\end{equation}
there exist an increasing sequence $\{s_k\}_{k\in\mathbb{N}}\subset(s_0,1)$ converging to $1$ and a function $u\in W^{1,p}_0(\Omega)$ such that
\[
\lim_{k\to\infty}\|u_{s_k}-u\|_{L^p(\Omega)}=0.
\]
\end{lemma}
\begin{proof}
By Poincar\'e inequality \eqref{spoin}, the estimate \eqref{ipotesi_utile} implies
\begin{equation}
\label{Lp}
\|u_s\|_{L^p(\Omega)}\le C_1,\qquad \mbox{ for every } s_0<s<1,
\end{equation}
for some $C_1=C_1(N,p,\mathrm{diam}(\Omega),L)>0$.
Moreover, again by \cite[Lemma A.1]{BraLinPar} and \eqref{ipotesi_utile} there exists a constant $C_2=C_2(N,p,L)>0$ such that
\begin{equation}
\label{power!}
\sup_{0<|\xi|<1}\int_{\mathbb{R}^N} \frac{|u_s(x+\xi)-u_s(x)|^p}{|\xi|^{s\,p}}\,dx\le C_2,\qquad \mbox{ for every } s_0<s<1.
\end{equation}
Since $s>s_0$, from the previous estimate, we can also infer
\begin{equation}
\label{s0}
\sup_{0<|\xi|<1}\int_{\mathbb{R}^N} \frac{|u_s(x+\xi)-u_s(x)|^p}{|\xi|^{s_0p}}\,dx\le C_2,\qquad \mbox{ for every } s_0<s<1.
\end{equation}
Estimates \eqref{Lp} and \eqref{s0} and the fact that $u_s\equiv0$ in $\mathbb{R}^N\setminus \Omega$ enables us to use the Riesz-Fr\'echet-Kolmogorov Compactness Theorem for $L^p$. Thus, there exists a sequence $\{u_{s_k}\}_{k\in\mathbb{N}}$ and $u\in L^p(\mathbb{R}^N)$ such that
\[
\lim_{k\to\infty} \|u_{s_k}-u\|_{L^p(\Omega)}=0.
\]
In order to conclude, we need to prove that $u\in W^{1,p}_0(\Omega)$. Up to a subsequence, we can suppose that $u_{s_k}$ converges almost everywhere. This implies that $u\equiv 0$ in $\mathbb{R}^N\setminus\Omega$. Moreover, thanks to Fatou Lemma we can pass to the limit in \eqref{power!} and obtain
\[
\sup_{0<|\xi|<1}\int_{\mathbb{R}^N} \frac{|u(x+\xi)-u(x)|^p}{|\xi|^{p}}\,dx\le C_2.
\]
This implies that the distributional gradient of $u$ is in $L^p(\mathbb{R}^N)$. Thus $u\in W^{1,p}(\mathbb{R}^N)$ and it vanishes almost everywhere in $\mathbb{R}^N\setminus\Omega$. Since $\Omega$ is Lipschitz, this finally implies that $u\in W^{1,p}_0(\Omega)$ by \cite[Propostion IX.18]{brezis}.
\end{proof}
The following result will complete the proof of Theorem \ref{thm:gamma}.
\begin{proposition}[$\Gamma-\liminf$ inequality]
\label{prop:liminf}
Given $\{s_k\}_{k\in\mathbb{N}}\subset\mathbb{R}$ an increasing sequence converging to $1$ and $\{u_k\}_{k\in\mathbb{N}}\subset L^p(\Omega)$ converging to $u$ in $L^p(\Omega)$, we have
\begin{equation}
\label{liminf}
\mathcal{E}_{1,p}(u)^p\le \liminf_{k\to \infty} \mathcal{E}_{s_k,p}(u_k)^p.
\end{equation}
\end{proposition}
\begin{proof}
The proof follows that of \cite[Lemma 7]{ADePM}. We start by observing that if
\[
\liminf_{k \to \infty}(1-s_k)\, [u_{k}]^p_{W^{{s_k},p}(\mathbb{R}^N)}= + \infty,
\] 
there is nothing to prove. Thus, let us suppose that 
\[
\liminf_{k \to \infty} (1-s_k)\, [u_{k}]^p_{W^{{s_k},p}(\mathbb{R}^N)}< + \infty,
\] 
this implies that for $k$ sufficiently large we have
\[
u_k\in \widetilde W^{s_k,p}_0(\Omega)\qquad \mbox{ and }\qquad (1-s_k)\, [u_k]^p_{W^{{s_k},p}(\mathbb{R}^N)}\le L,
\]
for some uniform constant $L>0$. By Lemma \ref{lm:utile}, we get that $u\in W^{1,p}_0(\Omega)$.
\vskip.2cm\noindent
We now continue the proof of \eqref{liminf}. For every measurable set $A \subset \Omega$ we define the absolutely continuous measure
\[ 
\mu(A):= \int_A |\nabla u|^p\,dy, 
\]
and we observe that, by Lebesgue's Theorem
\begin{equation}
\label{condizionelebesgue} 
\lim_{r \to 0^+} \frac{\mu(x+r\,Q)}{r^N} = \lim_{r\to 0^+}\frac{1}{r^N}\,\int_{x+r\,Q} |\nabla u|^p\,dy = |\nabla u(x)|^p \qquad \text{for a.\,e. }x \in \Omega,
\end{equation}
where as before $Q=(-1/2,1/2)^N$.
For a Borel set $E \subset \Omega$ we define
\[ 
\alpha_k(E) := (1-s_k)\,[u_k]^p_{W^{{s_k},p}(E)}\qquad \mbox{ and }\qquad 
\alpha(E) := \liminf_{k \to \infty} \alpha_k(E).
\]
For $x \in \Omega$, set $C_r(x):= x + r\,Q$. We claim that
\begin{equation} 
\label{stimaambrosio} 
\liminf_{r \to 0} \frac{\alpha(C_r(x))}{\mu(C_r(x))} \geq K(p,N), \qquad \text{for }\mu-\text{a.\,e. }x \in \Omega.
\end{equation}
In order to prove \eqref{stimaambrosio}, for every measurable function $v$, we introduce the notation
\[
v_{r,x}(y):=\frac{v(r\,y+x)-u(x)}{r},\qquad y\in Q.
\]
We keep on using the notation $\Psi_a(x)=\langle a,x\rangle$, for any given vector $a\in\mathbb{R}^N$.
Then we will prove \eqref{stimaambrosio} at any point $x \in \Omega$ such that 
\begin{equation}
\label{blow-up}
\lim_{r\searrow 0} \|u_{r,x}-\Psi_a\|_{L^{p}(Q)}=0,\qquad\mbox{ where } a = \nabla u(x),
\end{equation}
and such that \eqref{condizionelebesgue} holds. We recall that \eqref{blow-up} is true at almost every $x\in\Omega$ by \cite[Theorem 2, page 230]{EG}.
Therefore, to prove \eqref{stimaambrosio} it will be sufficient to show that
\[ 
\liminf_{r \to 0} \frac{\alpha(C_r(x))}{r^N} \geq K(p,N)\,|\nabla u(x)|^p.
\]
To this aim, let $r_j \searrow 0$ be a sequence such that
\[ \lim_{j \to \infty} \frac{\alpha(C_{r_j}(x))}{r_j^N} = \liminf_{r \to 0} \frac{\alpha(C_r(x))}{r^N}. \]
For any $j\in\mathbb{N}$ we can choose $k=k(j)$ so large that
\begin{enumerate}
\item[{\it i)}] $\alpha_{k(j)}(C_{r_j}(x)) \leq \alpha(C_{r_j}(x)) + r_j^{N+1}$;
\item[{\it ii)}] $\,r_j^{(1-s_{k(j)})\,p} \geq 1- 1/j$;
\item[{\it iii)}] $r_j^{-N-p}\, \|u_{k(j)}-u\|^p_{L^p(C_{r_j}(x))} <1/j$.
\end{enumerate}
Then, by using {\it i)}, the definitions of $\alpha_k$ and $(u_{k})_{r,x}$ and {\it ii)} we have
\begin{align*} 
\frac{\alpha(C_{r_j}(x))}{r_j^N} & \geq  \frac{\alpha_{k(j)}(C_{r_j}(x))}{r_j^N} - r_j\\
&= \frac{(1-s_{k(j)}) \,r_j^{N-s_{k(j)}\,p}\,r_{j}^p\,\Big[(u_{k(j)})_{r_j,x}\Big]^p_{W^{s_{k(j)},p}(Q)}}{r_j^N} - r_j \\ 
& \geq \left(1- \frac{1}{j}\right)\,(1-s_{k(j)})\,\Big[(u_{k(j)})_{r_j,x}\Big]^p_{W^{s_{k(j)},p}(Q)} - r_j. 
\end{align*}
On the other hand by {\it iii)} we have
\[ 
\Big\|(u_{k(j)})_{r_j,x}-u_{r_j,x}\Big\|_{L^p(Q)}<\left(\frac{1}{j}\right)^\frac{1}{p},
\]
while by \eqref{blow-up}
\[ 
\lim_{j\to\infty}\|u_{r_j,x} - \Psi_a\|_{L^p(Q)}= 0.
\]
Thus by triangle inequality we get that $(u_{k(j)})_{r_j,x}$ converges to $\Psi_a$ in $L^p(Q)$, with direction $a = \nabla u(x)$. This in turn implies
\[ 
\begin{split}
\lim_{j \to \infty} \frac{\alpha(C_{r_j}(x))}{r_j^N}&\geq \liminf_{j \to \infty} \left[\left(1- \frac{1}{j}\right)\,(1-s_{k(j)})\,\Big[(u_{k(j)})_{r_j,x}\Big]^p_{W^{s_{k(j)},p}(Q)} - r_j\right]\\
&\ge \Theta\left(p,N;\frac{\nabla u(x)}{|\nabla u(x)|}\right)\,|\nabla u(x)|^p\\
&=K(p,N)\,|\nabla u(x)|^p,\qquad \mbox{ for $\mu-$a.\,e. } x\in\Omega,
\end{split}
\]
thanks to the definition \eqref{gammapna} of $\Theta(p,N;a)$, property \eqref{nonunit} and Lemma \ref{lm:costanti}. This proves \eqref{stimaambrosio}. 
\vskip.2cm\par
The conclusion is exactly as in \cite[Lemma 7]{ADePM}. Let us consider for $\varepsilon > 0$ the following family of closed cubes
\[ 
\mathfrak{F} := \Big\{\overline{C_r(x)} \subset \Omega\,:\,(1+\varepsilon)\,\alpha\left(\overline{C_r(x)}\right) \geq K(p,N)\,\mu\left(\overline{C_r(x)}\right) \Big\}.
\]
By observing that
\[
\alpha\left(\overline{C_r(x)}\right)=\alpha(C_r(x))\qquad \mbox{ and }\qquad \mu\left(\overline{C_r(x)}\right)=\mu(C_r(x)),
\]
 and using \eqref{stimaambrosio}, we get that $\mathfrak{F}$ is a {\it fine Morse cover} (see \cite[Definition 1.142]{FL}) of $\mu-$almost all of $\Omega$, then we can apply a suitable version of Besicovitch Covering Theorem (see \cite[Corollary 1.149]{FL}) and extract a countable subfamily of disjoint cubes $\{C_i\}_{i \in I}\subset\mathfrak{F}$ such that $\mu(\Omega \setminus \cup_{i \in I} C_i)=0$. This yields
\begin{align*} 
K(p,N)\, \int_\Omega |\nabla u|^p\,dx & = K(p,N)\, \mu\left(\bigcup_{i \in I}C_i\right) = K(N,p) \sum_{i \in I} \mu(C_i)\\
&\leq (1+\varepsilon) \sum_{i \in I} \alpha(C_i) \leq (1+\varepsilon)\, \liminf_{k \to \infty} \sum_{i \in I} \alpha_k(C_i)\\
&\leq (1+\varepsilon)\, \liminf_{k \to \infty} 
 (1-s_k)\,[u_k]^p_{W^{{s_k},p}(\mathbb{R}^N)}.
\end{align*}
By the arbitrariness of $\varepsilon$ we get
\[
K(p,N) \int_\Omega |\nabla u|^p\,dx \leq \liminf_{k \to \infty} (1-s_k)\,[u_k]^p_{W^{{s_k},p}(\mathbb{R}^N)}.
\]
This concludes the proof.
\end{proof}

\subsection{A comment on dual norms}

By using Theorem \ref{thm:gamma}, we can prove a dual version of the Bourgain-Brezis-Mironescu result. The result of this section is not needed for the proof of Theorem \ref{main} and is placed here for completeness.
\begin{proposition}
Let $1<p<\infty$, for every $F\in L^{p'}(\Omega)$ we have 
\begin{equation}
\label{duali}
\lim_{s\nearrow 1} (1-s)^{-\frac{1}{p}}\,\|F\|_{W^{-s,p'}(\Omega)}=K(p,N)^{-\frac{1}{p}}\,\|F\|_{W^{-1,p'}(\Omega)}.
\end{equation}
\end{proposition}
\begin{proof}
We are going to use the variational characterization \eqref{duale} for dual norms.
By Corollary \ref{cor:gammacor}, the family of functionals
\begin{equation}
\label{modificati}
\mathcal{E}_{s,p}(u)^p-\int F\,u\,dx,\qquad u\in L^p(\Omega),
\end{equation}
$\Gamma-$converges to 
\[
\mathcal{E}_{1,p}(u)^p-\int_\Omega F\,u\,dx,\qquad u\in L^p(\Omega).
\]
We now observe that the functionals \eqref{modificati} are equi-coercive on $L^p(\Omega)$. Indeed, if $u\in L^p(\Omega)$ is such that
\begin{equation}
\label{equilo}
\mathcal{E}_{s,p}(u)^p-\int_\Omega F\,u\,dx\le M,
\end{equation}
this of course implies that $u\in \widetilde W^{s,p}_0(\Omega)$. Moreover, by Young inequality and \eqref{mandorlini} we have
\[
\begin{split}
M&\ge (1-s)\,[u]^p_{W^{s,p}(\mathbb{R}^N)}-\int_\Omega F\,u\,dx\\
&\ge (1-s)\,[u]^p_{W^{s,p}(\mathbb{R}^N)}-\frac{1}{p'}\,\left(\frac{\|F\|_{W^{-s,p}(\Omega)}}{(1-s)^\frac{1}{p}}\right)^{p'}-\frac{1-s}{p}\, [u]^p_{W^{s,p}(\mathbb{R}^N)}\\
&\ge \frac{1}{p'}\,(1-s)\,[u]^p_{W^{s,p}(\mathbb{R}^N)}-C\, \|F\|_{W^{-1,p'}(\Omega)},
\end{split}
\]
for a constant $C=C(N,p,\Omega)>0$ independent of $s$ (provided $s$ is sufficiently close to $1$). Thus from \eqref{equilo} we get
\[
(1-s)\,[u]^p_{W^{s,p}(\mathbb{R}^N)}\le M\,p'+C\,p'\, \|F\|_{W^{-1,p'}(\Omega)}.
\]
The desired equicoercivity in $L^p(\Omega)$ now follows from Lemma \ref{lm:utile}. In conclusion, from \eqref{duale}, the $\Gamma-$convergence and \eqref{duale1} we get
\[
\begin{split}
\lim_{s\nearrow 1}\frac{\|F\|_{W^{-s,p'}(\Omega)}}{(1-s)^\frac{1}{p}}&=\lim_{s\nearrow 1}\,p^\frac{1}{p}\,\left(-p'\min_{u\in L^p(\Omega)} \left\{\mathcal{E}_{s,p}(u)^p-\int_\Omega F\,u\, dx\right\}\right)^\frac{1}{p'}\\
&=p^\frac{1}{p}\,\left(-p'\min_{u\in L^p(\Omega)} \left\{\mathcal{E}_{1,p}(u)^p-\int_\Omega F\,u\,dx\right\}\right)^\frac{1}{p'}=\frac{\|F\|_{W^{-1,p}(\Omega)}}{K(p,N)^\frac{1}{p}},
\end{split}
\]
as desired.
\end{proof}
\begin{remark}
We recall the following dual characterization of $\|\cdot\|_{W^{-s,p'}(\Omega)}$ from \cite[Section 8]{BraLinPar}
\begin{equation}
\label{dualeBLP}
\|F\|_{W^{-s,p'}(\Omega)}=\min_{\varphi\in L^{p'}(\mathbb{R}^N\times\mathbb{R}^N)}\left\{\|\varphi\|_{L^{p'}(\mathbb{R}^N\times\mathbb{R}^N)}\, :\, R^*_{s,p}(\varphi)=F\mbox{ in }\Omega\right\},
\end{equation}
where $R^*_{s,p}$ is the adjoint of the linear continuous operator $R_{s,p}: \widetilde W^{s,p}_0(\Omega)\to L^p(\mathbb{R}^N\times \mathbb{R}^N)$  defined by
\[
R_{s,p}(u)(x,y)=\frac{u(x)-u(y)}{|x-y|^{\frac{N}{p}+s}},\qquad \mbox{ for every } u\in \widetilde W^{s,p}_0(\Omega).
\]
Formula \eqref{dualeBLP} is the nonlocal analog of the well-known duality formula
\[
\|F\|_{W^{-1,p'}(\Omega)}=\min_{V\in L^{p'}(\Omega;\mathbb{R}^N)}\left\{\|V\|_{L^{p'}(\Omega;\mathbb{R}^N)}\, :\, -\mathrm{div}\,V=F\mbox{ in }\Omega\right\}.
\]
\end{remark}
Then we end this section with the following curious convergence result.
\begin{corollary}
Let $1<p<\infty$ and $F\in L^{p'}(\Omega)$, then we have
\[
\begin{split}
\lim_{s\nearrow 1}\Big[(1-s)^{-\frac{1}{p}}&\min_{\varphi\in L^{p'}(\mathbb{R}^N\times\mathbb{R}^N)}\left\{\|\varphi\|_{L^{p'}(\mathbb{R}^N\times\mathbb{R}^N)}\, :\, R^*_{s,p}(\varphi)=F\mbox{ in }\Omega\right\}\Big]\\
&=K(p,N)^{-\frac{1}{p}}\, \min_{V\in L^{p'}(\Omega;\mathbb{R}^N)}\left\{\|V\|_{L^{p'}(\Omega;\mathbb{R}^N)}\, :\, -\mathrm{div}\,V=F\mbox{ in }\Omega\right\}.
\end{split}
\]
\end{corollary}

\section{Proof of Theorem~\ref{main}} 
\label{Proof}
\subsection{Convergence of the variational eigenvalues}

By Theorem \ref{thm:gamma}, we already know that
\[
\mathcal E_{1,p}(u)=\Big(\Gamma-\lim_{k\to\infty}\mathcal E_{s_k,p}\Big)(u),\qquad\mbox{ for all } u\in L^p(\Omega).
\]
For every $1<p<\infty$, let us define the functional $g_{p}:L^p(\Omega)\to[0,\infty)$ by
$$
g_{p}(u):=\|u\|_{L^p(\Omega)}.
$$  
We now observe that for every increasing sequence $k_n$ and for 
any sequence $\{u_n\}_{n\in\mathbb{N}}\subset L^p(\Omega)$ such that
\[
M:=\sup_{n\in\mathbb N}\mathcal E_{s_{k_n},p}(u_n)<+\infty,
\]
there exists a subsequence $\{u_{n_j}\}_{j\in\mathbb{N}}$ such that
\[
\lim_{j\to\infty}g_{p}(u_{n_j})=g_{p}(u).
\]
Indeed, this is a consequence of Lemma \ref{lm:utile}.
Then the functionals $\mathcal E_{s,p}$ and $g_{p}$ satisfy all the assumptions in \cite[Corollary 4.4]{DM14}, which implies
\begin{equation}
\label{gamma-conv}
\lim_{k\to\infty}\left(\inf_{K\in\mathcal K_{m,p}(\Omega)}\sup_{u\in K}\mathcal E_{s_k,p}(u)\right)=\inf_{K\in\mathcal K_{m,p}(\Omega)}\sup_{u\in K}\mathcal E_{1,p}(u),
\end{equation}
where
\[
\mathcal{K}_{m,p}(\Omega)=\Big\{K\subset \{u\,:\, g_{p}(u)=1\}\, :\, K \mbox{ compact and symmetric}, \, i(K)\ge m\Big\}.
\]
In order to conclude, we only need to show that the minimax values with respect to the $\widetilde W^{s,p}_0(\Omega)-$topology are equal to those with respect to the weaker topology $L^p(\Omega)$. Observe now that, for every $b\in\mathbb{R}$, the
restriction of $g_p$ to $\{u\in L^p(\Omega): {\mathcal  E}_{s,p}(u)\leq b\}$ is continuous: for $s=1$ this is classical, while for $0<s<1$ we can appeal for example to \cite[Theorem 2.7]{BraLinPar}.
Whence, \cite[Corollary 3.3]{DM14} yields
\begin{align}
\label{id1}
\inf_{K\in\mathcal K_{m,p}(\Omega)}\sup_{u\in K}\mathcal E_{1,p}(u) &=\inf_{K\in\mathcal W^{1}_{m,p}(\Omega)}\sup_{u\in K}\mathcal E_{1,p}(u), \\
\label{id2}
\inf_{K\in\mathcal K_{m,p}(\Omega)}\sup_{u\in K}\mathcal E_{s_k,p}(u)&=\inf_{K\in\mathcal W^{s_k}_{m,p}(\Omega)}\sup_{u\in K}\mathcal E_{s_k,p}(u),
\end{align}
where we recall that $\mathcal{W}^s_{m,p}(\Omega)$ has been defined in \eqref{doppiwu}.
By using \eqref{id1} and \eqref{id2} in \eqref{gamma-conv}, 
the assertion follows by definition of $\lambda^{s}_{m,p}(\Omega)$ and $\lambda^1_{m,p}(\Omega)$. 

\subsection{Convergence of the eigenfunctions}

For every $s\in(0,1)$, let $u_s\in \widetilde W^{s,p}_0(\Omega)$ be an eigenfunction corresponding to the variational eigenvalue $\lambda^s_{m,p}(\Omega)$, normalized by $\|u_s\|_{L^p(\Omega)}=1$. Then it verifies
\[
(1-s)\,[u_s]^p_{W^{s,p}(\mathbb{R}^N)}=(1-s)\,\lambda^s_{m,p}(\Omega).
\] 
The convergence of the eigenvalues, which has been proved in the previous subsection, implies that 
\begin{equation}
\label{carlona}
(1-s)\,[u_s]^p_{W^{s,p}(\mathbb{R}^N)}\le K(p,N)\,\lambda^1_{m,p}(\Omega)+1,
\end{equation}
up to choosing $1-s$ sufficiently small.
By appealing again to Lemma \ref{lm:utile}, this in turn implies that there exists a sequence $\{s_k\}_{k\in\mathbb{N}}$ with $s_k\nearrow 1$ such that the corresponding sequence of eigenfunctions $\{u_{s_k}\}_{k\in\mathbb{N}}$ converges strongly in $L^p$ to a function $u\in W^{1,p}_0(\Omega)$. By strong convergence, we still have $\|u\|_{L^p(\Omega)}=1$. 
\par
In order to prove that $u$ is an eigenfunction of the local problem, let us notice
that each $u_{s_k}$ weakly solves
\[
(-\Delta_p)^{s_k} u=\lambda^{s_k}_{m,p}(\Omega)\,|u_{s_k}|^{p-2}\,u_{s_k},\quad \mbox{ in }\Omega,\qquad u=0,\quad \mbox{ in }\mathbb{R}^N\setminus \Omega.
\]
Thus it is the unique minimizer of the following strictly convex problem
\[
\min_{v\in L^p(\Omega)} \left\{\mathcal{E}_{s_k,p}(v)^p+p\,\int_\Omega F_{s_k}\, v\,dx\right\},
\]
where 
\[
F_{s_k}=-(1-s_k)\, \lambda^{s_k}_{m,p}(\Omega)\,|u_{s_k}|^{p-2}\,u_{s_k}\in L^{p'}(\Omega).
\]
Observe that the sequence $\{F_{s_k}\}_{k\in\mathbb{N}}$ converges strongly in $L^{p'}(\Omega)$ to the function
\begin{equation}
\label{F}
F=-K(p,N)\, \lambda^1_{m,p}(\Omega)\,|u|^{p-2}\,u,
\end{equation}
thanks to the strong convergence of $\{u_{s_k}\}_{k\in\mathbb{N}}$ and to the first part of the proof. By appealing to the $\Gamma-$convergence result of Corollary \ref{cor:gammacor}, we thus get that $u$ is a solution (indeed the unique, again by strict convexity) of the limit problem
\[
\min_{v\in L^p(\Omega)} \left\{\mathcal{E}_{1,p}(v)^p+p\,\int_\Omega F\, v\,dx\right\},
\]
with $F\in L^{p'}(\Omega)$ defined in \eqref{F}. As a solution of this problem, $u$ has to satisfy the relevant Euler-Lagrange equation, i.e. $u$ weakly solves
\[
-\Delta_p u=\lambda^{1}_{m,p}(\Omega)\,|u|^{p-2}\,u,\quad \mbox{ in }\Omega,\qquad u=0,\quad \mbox{ on }\partial\Omega.
\]
This proves that the renormalized eigenfunctions $\{u_{s_k}\}_{k\in\mathbb{N}}$ converges strongly in $L^{p}(\Omega)$ to an eigenfunction $u$ corresponding to $\lambda^1_{m,p}(\Omega)$ having unit norm.
\par
In order to improve the convergence in $\widetilde W^{t,q}_0(\Omega)$ for every $p\le q<\infty$ and every $t<p/q$, it is now sufficient to use the interpolation inequality of Proposition \ref{prop:interpolation} with $r=+\infty$ so that $\alpha=s_k\,p/q$. Observe that since $s_k$ is converging to $1$, if $t<p/q$ we can always suppose that
\[
t<s_k\,\frac{p}{q},
\]
up to choosing $k$ large enough.
This yields for a constant\footnote{We may notice that the constant $C$ degenerates as $t$ approaches $p/q$.} $C=C(N,p,q,t)>0$ which varies from a line to another
\[
\begin{split}
\lim_{k\to\infty} &t^\frac{1}{q}\,[u_{s_k}-u]_{W^{t,q}(\mathbb{R}^N)}\\
&\le C\,\lim_{k\to\infty}\|u_{s_k}-u\|_{L^p(\Omega)}^{\frac{p}{q}\,\left(1-\frac{q}{p}\,\frac{t}{s_k}\right)}\,\|u_{s_k}-u\|^{\left(1-\frac{q}{p}\right)}_{L^\infty(\Omega)}\\
&\times\left((1-s_k)^\frac{1}{p}\,[u_{s_k}-u]_{W^{s_k,p}(\mathbb{R}^N)}\right)^{\frac{t}{s_k}\,}\\
&\le C\,\lim_{k\to\infty}\|u_{s_k}-u\|_{L^p(\mathbb{R}^N)}^{\frac{p}{q}\,\left(1-\frac{q}{p}\,\frac{t}{s_k}\right)}\\
&\times\left(\|u_{s_k}\|_{L^\infty(\Omega)}+\|u\|_{L^\infty(\Omega)}\right)^{\left(1-\frac{q}{p}\right)}\\
&\times\left((1-s_k)^\frac{1}{p}\,[u_{s_k}]_{W^{s_k,p}(\mathbb{R}^N)}+(1-s_k)\,[u]_{W^{s_k,p}(\mathbb{R}^N)}\right)^\frac{t}{s_k}.
\end{split}
\]
If we use \eqref{boundep} (if $1<p< N$), \eqref{cresima} (if $p=N$) or \eqref{bounder} (if $p>N$) to bound the $L^\infty$ norms, \eqref{carlona} and Proposition \ref{norm-cont} to bound the Gagliardo semi-norms, we finally get
\[
\lim_{k\to\infty}t^\frac{1}{q}\,[u_{s_k}-u]_{W^{t,q}(\mathbb{R}^N)}\le  C\, \lim_{k\to\infty}\|u_{s_k}-u\|_{L^p(\mathbb{R}^N)}^{\frac{p}{q}\,\left(1-\frac{q}{p}\,\frac{t}{s_k}\right)}=0,
\]
as desired.
\begin{remark}[Pushing the convergence further]
In the previous result, we used that the initial convergence in $L^p$ norm can be ``boosted'' by combining suitable interpolation inequalities and regularity estimates exhibiting the correct scaling in $s$. Thus, should one obtain that eigenfunctions are more regular with good a priori estimates, the previous convergence result could still be improved. Though it is known that eigenfunctions are continuous for every $1<p<\infty$ and $0<s<1$ (see \cite{Kuusi,IaMosSq}), unfortunately the above mentioned results do not provide estimates with an explicit dependence on $s$ and thus we can not directly use them. 
\par
In the case $p=2$, regularity estimates of this type can be found in \cite[Lemma 4.4]{CS} for bounded solutions of the equation in the whole space
\[
(-\Delta)^s u=f(u)\qquad \mbox{ in }\mathbb{R}^N,
\]
where $f$ is a (smooth) nonlinearity. For such an equation, the authors prove Schauder-type estimates for the solutions, with constants independent of $s$ (provided $s>s_0>0$).
\end{remark}

\appendix

\section{Courant vs.\ Ljusternik-Schnirelmann}
\label{app:cls}

Here we prove that for $p=2$ the variational eigenvalues defined by the Ljusternik-Schnirelman procedure \eqref{lambdas} coincide with the usual eigenvalues coming from Spectral Theory (see Theorem \ref{thm:uguali!} below). Thus in particular for $p=2$ definition \eqref{lambdas} give {\it all} the eigenvalues.
This fact seems to belong to the folklore of Nonlinear Analysis, but since we have not been able to find a reference in the literature, we decided to include this Appendix.
\vskip.2cm
Let $H_1\subset H_2$ be two separable infinite dimensional Hilbert spaces, endowed with scalar products $\langle \cdot,\cdot\rangle_{H_i}$ and norms
\[
\|u\|_{H_i}=\sqrt{\langle u,u\rangle_{H_i}},\qquad u\in H_i,\ i=1,2.
\]
On the space $H_1$ is defined a symmetric bilinear form $\mathcal{Q}:H_1\times H_1\to\mathbb{R}$.
We assume the following:
\begin{itemize}
\item[1.] the inclusion $\mathcal{I}:H_1\to H_2$ is a continuous and compact linear operator;
\vskip.2cm 
\item[2.] $\mathcal{Q}$ is continuous and coercive, i.e. for $C\ge 1$
\[
\frac{1}{C}\,\|u\|^2_{H_1}\le\mathcal{Q}[u,u],\qquad \mathcal{Q}[u,v]\le C\,\|u\|_{H_1}\,\|v\|_{H_1}\qquad u,v\in H_1.
\]
Thus $\mathcal{Q}$ defines a scalar product on $H_1$, whose associated norm is equivalent to $\|\cdot\|_{H_1}$.
\end{itemize}
We set
\[
\mathcal{S}=\{u\in H_1\, :\, \|u\|_{H_2}=1 \}.
\]
Then the restriction of the functional $u\mapsto \mathcal{Q}[u,u]$ to $\mathcal{S}$ has countably many critical values $0<\lambda_1\le \lambda_2\le\dots\le \lambda_m\le \dots \nearrow +\infty$, with associated a sequence of critical points $\{\varphi_n\}_{n\in\mathbb{N}}\subset \mathcal{S}$ defining a Hilbertian basis of $H_2$.\footnote{These are indeed the inverses of the eigenvalues of the {\it resolvent operator} $\mathcal{R}:H_2\to H_2$ defined by:
\[
\mbox{ for } f\in H_2,\quad \mathcal{R}(f)\in H_1\subset H_2 \mbox{ is the unique solution of }
\]
\[
\mathcal{Q}[\mathcal{R}(f),u]=\langle f,u\rangle_{H_2},\qquad \mbox{ for every }u \in H_1.
\]
The hypotheses above guarantee that $\mathcal{R}$ is a well-defined compact, positive and self-adjoint linear operator. Then discreteness of the spectrum follows from the Spectral Theorem, see for example \cite[Theorem 1.2.1]{He}.}   The critical point $\varphi_i$ satisfies
\[
\mathcal{Q}[\varphi_i,u]=\lambda_i\, \langle \varphi_i,u\rangle_{H_2},\qquad \mbox{ for every }\varphi\in H_1,
\]
so in particular
\[
\mathcal{Q}[\varphi_i,\varphi_j]=\lambda_i\,\delta_{ij},\qquad i,j\in\mathbb{N}\setminus\{0\}.
\]
These critical points have a variational characterization: indeed, if we introduce
for every $m\in\mathbb{N}\setminus\{0\}$
\[
\mathcal{E}_{m}=\{E\subset H_1\, :\, E \mbox{ vector space with } \mathrm{dim}(E)\ge m\},
\]
and
\[
\mathcal{F}_{m}=\{F\subset \mathcal{S}\, :\, F=E\cap \mathcal{S}\mbox{ for some } E\in \mathcal{E}_{m}\},
\]
then we have
\[
\lambda_m=\min_{F\in\mathcal{F}_{m}} \max_{u\in F} \mathcal{Q}[u,u].
\]
A minimizer for the previous problem is given by
\begin{equation}
\label{Fm}
F_m:=\mathrm{Span}\big\{\varphi_1,\dots,\varphi_m\big\}.
\end{equation}
We also recall that the eigenvalues can be also characterized as
\[
\lambda_m=\min_{u\in H_1\setminus\{0\}}\left\{\frac{\mathcal{Q}[u,u]}{\|u\|^2_{H_2}}\, :\, \langle u,\varphi_i\rangle_{H_2}=0,\ i=1,\dots,m-1\right\}.
\]
\begin{lemma}
\label{lm:codroipo}
Let $m\in\mathbb{N}\setminus\{0\}$, we define 
\[
\mathcal{W}_{m}=\{K\subset \mathcal{S}\, :\, K \mbox{ compact and symmetric},\, i(K)\ge m\},
\]
where $i$ is the {\it Krasnosel'ski\u{\i} genus}, see definition \eqref{krasno}.
Then for every $K\in\mathcal{W}_{m}$ we have
\begin{equation}
\label{complement}
K\cap F_{m-1}^\bot\not=\emptyset.
\end{equation}
Here $F_{m-1}$ is defined as in \eqref{Fm} and orthogonality is intended in the $H_2$ sense.
\end{lemma}
\begin{proof}
We proceed by contradiction. Let us assume that there exists $K\in\mathcal{W}_{m}$ such that \eqref{complement} is not true, this implies that
\begin{equation}
\label{diversi}
\sum_{j=1}^{m-1} \langle u,\varphi_j\rangle_{H_2}^2\not=0,\qquad \mbox{for every } u\in K.
\end{equation}
We can now define the following map $\Phi:K\to \mathbb{S}^{m-2}$ by
\[
\begin{array}{ccccc}
& & u& \mapsto & \displaystyle\sum_{j=1}^{m-1} \mathbf{e}_j\, \langle u,\varphi_j\rangle_{H_2}\,\left[\displaystyle \sum_{j=1}^{m-1} \langle u,\varphi_j\rangle_{H_2}^2\right]^{-\frac{1}{2}},
\end{array}
\]
where $\mathbf{e}_j$ is the $j-$th versor of the canonical basis. Thanks to \eqref{diversi}, the previous map is well-defined, continuous and odd. This contradicts the fact that $K$ has genus greater or equal than $m$ and thus \eqref{complement} holds true.
\end{proof}
\begin{theorem}
\label{thm:uguali!}
For every $m\in\mathbb{N}\setminus\{0\}$, we have
\begin{equation}
\label{uguali!}
\lambda_m=\inf_{K\in\mathcal{W}_{m}} \max_{u\in K} \mathcal{Q}[u,u].
\end{equation}
\end{theorem}
\begin{proof}
The inequality 
\[
\lambda_m\ge \inf_{K\in\mathcal{W}_{m}} \max_{u\in K} \mathcal{Q}[u,u],
\]
easily follows the fact that $\mathcal{F}_{m}\subset\mathcal{W}_{m}$ (see \cite[Chapter 2, Proposition 5.2]{St}).
\vskip.2cm\noindent
In order to prove the reverse inequality, for every $\varepsilon>0$ let $K_\varepsilon\in\mathcal{W}_{m}$ be such that
\[
\max_{u\in K_\varepsilon} \mathcal{Q}[u,u]\le\inf_{K\in\mathcal{W}_{m}} \max_{u\in K} \mathcal{Q}[u,u]+\varepsilon.
\]
By Lemma \ref{lm:codroipo}, $K_\varepsilon$ is such that $K_\varepsilon \cap F_{m-1}^\bot\not=\emptyset$.
In particular, there exists $v\in K_\varepsilon$ such that
\[
v=\sum_{j=m}^\infty \alpha_{j}\,\varphi_j,\qquad \mbox{ with } \sum_{j=m}^\infty \alpha_j^2=1.
\]
We have
\[
\max_{u\in K_\varepsilon}\mathcal{Q}[u,u]\ge \mathcal{Q}[v,v]=\sum_{j=m}^\infty \alpha_j^2\, \lambda_j\ge \lambda_m.
\]
This in turn implies
\[
\lambda_m\le\inf_{K\in\mathcal{W}_{m}} \max_{u\in K}\mathcal{Q}[u,u]+\varepsilon,
\]
and by the arbitrariness of $\varepsilon$ we get the conclusion.
\end{proof}

\section{A density result}

For completeness, we present the density result below. This permits to infer that the space
\[
\Big\{u:\mathbb{R}^N\to\mathbb{R}\, :\, [u]_{W^{s,p}(\mathbb{R}^N)}<+\infty\,\mbox{ and }\, u=0 \mbox{ in } \mathbb{R}^N \setminus \Omega\Big\}
\]
coincides with the completion of $C^\infty_0(\Omega)$ with respect to $[\,\cdot\,]_{W^{s,p}(\mathbb{R}^N)}$. For a fairly more general result obtained by means of a different proof, we also refer to \cite[Theorem 6]{FSV}.
\begin{proposition}
\label{prop:density}
Let $1<p<\infty$ and $s\in(0,1)$. Let $\Omega\subset\mathbb{R}^N$ be an open bounded set with Lipschitz boundary. For every measurable function $u:\mathbb{R}^N\to\mathbb{R}$ such that
\[
[u]_{W^{s,p}(\mathbb{R}^N)}<+\infty\qquad \mbox{ and }\qquad u= 0\mbox{ in }\mathbb{R}^N\setminus\Omega,
\]
 there exists a sequence $\{\varphi_n\}_{n\in\mathbb{N}}\subset C^\infty_0(\Omega)$ such that
\begin{equation}
\label{LS}
\lim_{n\to\infty} [\varphi_n-u]_{W^{s,p}(\mathbb{R}^N)}=0.
\end{equation}
\end{proposition}
\begin{proof}
The proof is similar to that of \cite[Lemma 2.3]{BraLinPar}, concerning the case $p=1$. By using \cite[Lemma 3.2]{littig_schuricht2014}, the regularity of $\Omega$ implies that there exists a family of diffeomorphisms $\Phi_\varepsilon:\mathbb{R}^N\to\mathbb{R}^N$ with inverses $\Psi_\varepsilon$ such that:
\begin{itemize}
\item we have
\[
\lim_{\varepsilon\to 0^+} \|D\Phi_\varepsilon-\mathrm{Id}\|_{L^\infty}+\|\Phi_\varepsilon-\mathrm{Id}\|_{L^\infty}=0,
\]
and 
\[
\lim_{\varepsilon\to 0^+} \|D\Psi_\varepsilon-\mathrm{Id}\|_{L^\infty}+\|\Psi_\varepsilon-\mathrm{Id}\|_{L^\infty}=0;
\]
\vskip.2cm
\item $\Omega_\varepsilon:=\Phi_\varepsilon(\overline\Omega)\Subset\Omega$ for all $\varepsilon\ll1$.
\end{itemize}
We then define the sequence $\varphi_n=(u\circ\Psi_{1/n})\ast\varrho_n$, where $\varrho_n$ is a positive convolution kernel such that $\|\varrho_n\|_{L^1}=1$ and chosen so that $\varphi_n$ has compact support in $\Omega$. Then by construction $\varphi_n\in C^\infty_0(\Omega)$ and
\[
\lim_{n\to\infty} \|\varphi_n-u\|_{L^p(\mathbb{R}^N)}=0.
\]  
By Fatou Lemma, this also implies that
\[
\liminf_{n\to\infty}\, [\varphi_n]_{W^{s,p}(\mathbb{R}^N)}\ge [u]_{W^{s,p}(\mathbb{R}^N)}.
\]
Moreover, we have
\begin{equation}
\label{jensen}
[\varphi_n]^p_{W^{s,p}(\mathbb{R}^N)}=[(u\circ\Psi_{1/n})\ast\varrho_n]^p_{W^{s,p}(\mathbb{R}^N)}\le [u\circ\Psi_{1/n}]^p_{W^{s,p}(\mathbb{R}^N)},
\end{equation}
which follows by using convexity of $\tau\mapsto \tau^p$ and Jensen's inequality with respect to the measure $\varrho_n\,dx$. Then we use that
\[
[u\circ\Psi_{1/n}]^p_{W^{s,p}(\mathbb{R}^N)}=\int_{\mathbb{R}^N}\int_{\mathbb{R}^N} \frac{|u(z)-u(w)|^p\,|J\Phi_{1/n}(z)|\,|J\Phi_{1/n}(w)|}{|\Phi_{1/n}(z)-\Phi_{1/n}(w)|^{N+s\,p}}\, dz\, dw, 
\]
which follows by a simple change of variables $(z,w)=(\Psi_{1/n}(x),\Psi_{1/n}(y))$, where $J\Phi_{1/n}$ denotes the Jacobian determinant. Observe that by construction
\[
|\Phi_{1/n}(z)-\Phi_{1/n}(w)|\ge M_1\, |z-w|\qquad \mbox{ and }\qquad |J\Phi_{1/n}(z)|\le M_2,
\]
for some $M_1>0$ and $M_2\ge 1$ independent of $n$. Thus we can apply Lebesgue Dominated Convergence Theorem and keeping into account \eqref{jensen}, we can infer that
\[
\limsup_{n\to\infty}\, [\varphi_n]_{W^{s,p}(\mathbb{R}^N)}\le [u]_{W^{s,p}(\mathbb{R}^N)},
\]
as well. In conclusion, we get that
\begin{equation}
\label{normeanorme}
\lim_{n\to\infty}\, [\varphi_n]_{W^{s,p}(\mathbb{R}^N)}= [u]_{W^{s,p}(\mathbb{R}^N)}.
\end{equation}
In order to conclude, by \eqref{normeanorme} the sequence
\[
\phi_n(x,y):=\frac{\varphi_n(x)-\varphi_n(y)}{|x-y|^{\frac{N}{p}+s}},
\]
is bounded in $L^p(\mathbb{R}^N\times\mathbb{R}^N)$ and it weakly converges to
\[
\phi(x,y):=\frac{u(x)-u(y)}{|x-y|^{\frac{N}{p}+s}}\in L^p(\mathbb{R}^N\times\mathbb{R}^N).
\]
By using this, \eqref{normeanorme} and uniform convexity of the $L^p$ norm, we get the desired result.
\end{proof}


\section*{Acknowledgments} 
We thank Guido De Philippis and Sunra Mosconi for some useful bibliographical comments and Matias Reivesle for a discussion on Poincar\'e inequalities.
Part of this paper was written during a visit of L.\,B. and E.\,P. in Verona and of M.\,S. in Marseille, as well as during the XXV Italian Workshop on Calculus of Variations held in Levico Terme, in February 2015. Hosting institutions and organizers are gratefully acknowledged.


\end{document}